\documentclass[reqno,11pt]{amsart}

\usepackage{graphicx}         
\usepackage{amsfonts}
\usepackage{amssymb}
\usepackage[latin1]{inputenc}
\usepackage[all]{xy}
\usepackage{tikz-cd}        

\newtheorem{thm}{Theorem}[section]

\newtheorem{prop}[thm]{Proposition}
\theoremstyle{definition}
\newtheorem{defn}[thm]{Definition}
\theoremstyle{remark}

\newtheorem{remark}[thm]{Remark}
\newtheorem{example}[thm]{Example}

\renewcommand{\a}{\alpha}
\renewcommand{\b}{\beta}

\newcommand{\s}{\sigma}

\newcommand{\PB}{\left\{\cdot\,,\cdot\right\}}
\newcommand{\pB}[1]{\left\{#1,\cdot\right\}}

\newcommand{\pb}[1]{\left\{#1\right\}}

\newcommand{\lb}[1]{\[#1\]}
\newcommand{\LB}{[\cdot\,,\cdot]}

\renewcommand{\(}{\left(}
\renewcommand{\)}{\right)}
\renewcommand{\[}{\left[}
\renewcommand{\]}{\right]}

\newcommand{\cA}{\mathcal A}
\newcommand{\cB}{\mathcal B}

\newcommand{\X}{ X}
\newcommand{\fg}{\mathfrak g}

\newcommand{\bbQ}{\mathbb Q}
\newcommand{\bbC}{\mathbb C}
\newcommand{\bbF}{\mathbb F}

\newcommand{\bbN}{\mathbb N}

\newcommand{\ux}{\underline{x}}
\newcommand{\uy}{\underline{y}}

\newcommand{\Gr}{\mathop{\rm gr}\nolimits}
\newcommand{\Sym}{\mathop{\rm Sym}\nolimits}
\newcommand{\omiet}[2]{\check{\underline{#1}}^{#2}}

\renewcommand{\Im}{\mathop{\rm Im}\nolimits}
\newcommand{\leqs}{\leqslant}
\newcommand{\geqs}{\geqslant}

\newcommand{\pp}[2]{\frac{\partial#1}{\partial#2}}

\newcommand{\diff}{{\rm d }}

\newcommand{\End}{{\mathop{\rm End}}}

\newcommand{\Id}{\mathrm{Id}}
\DeclareMathOperator{\PBW}{PBW}

\newif\ifprivate
\privatefalse

 \numberwithin{equation}{section}

\def\???{\ifprivate {\bf {???}} \marginpar{{\Huge {\bf ?}}}\else \fi}
\numberwithin{equation}{section}

\begin{document}  

\nocite{*} 

\title[Poisson enveloping algebras and PBW]{Poisson enveloping algebras and the Poincaré-Birkhoff-Witt theorem}

\author[Lambre]{Thierry Lambre} \address{Thierry Lambre, Laboratoire de Mathématiques Blaise Pascal, UMR 6620
  Université Clermont Auvergne, CNRS, Campus des Cézeaux, 3 place Vasarely, 63178 Aubière cedex,
  France}\email{thierry.lambre@uca.fr}

\author[Ospel]{Cyrille Ospel${}^\dagger$} \address{Cyrille Ospel, LaSIE, UMR 7356 CNRS, Université de La Rochelle,
  Av. M. Crépeau, 17042 La Rochelle cedex 1, La Rochelle, France}\email{cospel@univ-lr.fr}

\author[Vanhaecke]{Pol Vanhaecke${}^\dagger$} \address{Pol Vanhaecke, Laboratoire de Math\'ematiques et
  Applications, UMR 7348 CNRS, Universit\'e de Poitiers, 11 Boulevard Marie et Pierre Curie, Téléport 2 - BP 30179,
  86 962 Chasseneuil Futuroscope Cedex, France}\email{pol.vanhaecke@math.univ-poitiers.fr}

\thanks{${}^\dagger$Partially supported by the Federation MIRES FR CNRS 3423}

\date{\today}


\begin{abstract}
Poisson algebras are, just like Lie algebras, particular cases of Lie-Rinehart algebras. The latter were introduced
by Rinehart in his seminal 1963 paper, where he also introduces the notion of an enveloping algebra and proves ---
under some mild conditions --- that the enveloping algebra of a Lie-Rinehart algebra satisfies a
Poincaré-Birkhoff-Witt theorem (PBW theorem). In the case of a Poisson algebra $(\cA,\cdot,\PB)$ over a commutative
ring $R$ (with unit), Rinehart's result boils down to the statement that if $\cA$ is \emph{smooth} (as an algebra),
then $\Gr (U(\cA))$ and $\Sym(\Omega(\cA))$ are isomorphic as graded algebras; in this formula, $U(\cA)$ stands for
the Poisson enveloping algebra of $\cA$ and $\Omega(\cA)$ is the $\cA$-module of Kähler differentials of $\cA$
(viewing $\cA$ as an $R$-algebra). In this paper, we give several new constructions of the Poisson enveloping
algebra in some general and in some particular contexts. Moreover, we show that for an important class of
\emph{singular} Poisson algebras, the PBW theorem still holds. In geometrical terms, these Poisson algebras
correspond to (singular) Poisson hypersurfaces of arbitrary smooth affine Poisson varieties. Throughout the paper
we give several examples and present some first applications of the main theorem; applications to deformation
theory and to Poisson and Hochschild (co-) homology will be worked out in a future publication.
\end{abstract}

\maketitle

\setcounter{tocdepth}{2}

\tableofcontents

\section{Introduction}

Poisson brackets first appeared in classical mechanics as a tool for constructing new constants of motion from
given ones. Their importance was soon emphasized by the discovery that Poisson-commutativity is the key ingredient
of (Liouville) integrability. Since then, they have also played a major role in quantization theories (geometric
quantization, deformation quantization, quantum groups, \dots) and in many other parts of mathematical physics,
such as string theory. Formalizing the properties of the Poisson bracket leads on the geometrical side to the
definition of a Poisson manifold as being a manifold equipped with a bivector field whose induced operation on
smooth functions (the Poisson bracket!) is a Lie bracket. In purely algebraic terms, a Poisson algebra
$(\cA,\cdot,\PB)$ over some ring $R$ (always assumed commutative and unitary) is an $R$-module, equipped with a
commutative unitary algebra structure $''\cdot''$ and with a Lie algebra structure $\PB$, satisfying the following
compatibility relation, valid for all $a_1,a_2,a_3\in\cA$:
\begin{equation*}
       \pb{a_1\cdot a_2,a_3}=a_1\cdot \pb{a_2,a_3}+a_2\cdot \pb {a_1,a_3}\;.
\end{equation*}%
The algebra of smooth functions on a Poisson manifold is an important class of Poisson algebras, but even in the
realm of geometry, in particular in Lie theory and in algebraic geometry, one is soon led to considering singular
varieties, which are equipped with a Poisson structure. The simplest example of a singular Poisson variety is the
cone $X_1^2+X_2^2+X_3^2=0$ in $\bbC^3$, whose algebra of regular functions
\begin{equation*}
  \cB:=\frac{\bbC[X_1,X_2,X_3]}{(X_1^2+X_2^2+X_3^2)}
\end{equation*}%
is equipped with the Poisson structure, defined by the following three Poisson brackets:
\begin{equation}\label{eq:pb_intro}
  \pb{X_1,X_2}=X_3,\quad  \pb{X_2,X_3}=X_1,\quad  \pb{X_3,X_1}=X_2\;.
\end{equation}%
There is a natural notion of a Poisson module over a Poisson algebra, and hence of a Poisson enveloping algebra of
a Poisson algebra; the latter allows one to view a Poisson module over a Poisson algebra $\cA$ as a module over the
Poisson enveloping algebra $U(\cA)$ of $\cA$. Several existence proofs of the Poisson enveloping algebra are known,
see for example \cite{rinehart,huebsch,MR1683858}. For the case of (possibly modified) Lie-Poisson algebras we will
give two new constructions in Sections \ref{par:envel_lie-poisson} and \ref{par:envel_lie-poisson_2}. We also give
a new description of the Poisson enveloping algebra $U(\cB)$ of a quotient Poisson algebra $\cB=\cA/I$ in terms of
$U(\cA)$. This construction, which is valid for an arbitrary Poisson ideal $I$ in an arbitrary Poisson algebra
$\cA$, has some similarities with smash products, which are known to provide a construction of the Poisson
enveloping algebra of any Poisson algebra (see Section \ref{general_Uoisson}), but is yet quite different. This
description will play in important role in Section \ref{sec:PBW}.

A classical important result about the enveloping algebra of a Lie algebra is the Poincaré-Birkhoff-Witt theorem,
which states that when $\fg$ is a Lie algebra over a field of characteristic zero, there is a natural graded
algebra isomorphism between $\Sym\fg$ and the graded algebra of $U_{Lie}(\fg)$. For a \emph{smooth} Poisson algebra
$\cA$ a similar result holds: denoting by $\Omega(\cA)$ the $\cA$-module of Kähler differentials on $\cA$, there is
a natural graded algebra isomorphism between $\Sym_\cA\Omega(\cA)$ and the graded algebra of $U(\cA)$. The proof of
this very general result can be traced back to Rinehart \cite{rinehart}, which proves the PBW theorem for what we
call now \emph{Lie-Rinehart algebras}, of which Poisson algebras can be viewed as a special case (the latter fact
is an observation by Weinstein \cite{Weinstein,CD_Weinstein}. To be precise, when $\cA$ is a Poisson algebra, the
pair $(\cA,\Omega(\cA))$ is a Lie-Rinehart algebra; also, saying that $\cA$ is a \emph{smooth} algebra means in
algebraic terms that $\Omega(\cA)$ is a projective $\cA$-module, while it corresponds geometrically (for fields
such as $\bbC$) to the absence of singular points in the variety whose algebra of regular functions is $\cA$.

We show in this paper that the PBW theorem still holds for a large class of \emph{singular} (non-smooth) Poisson
algebras, including the algebra of functions of any irreducible Poisson hypersurface (possibly singular) of an
arbitrary smooth Poisson variety.


\begin{thm}\label{thm:intro_main}
  Suppose that the $\cA$ is a smooth Poisson algebra and that $I$ is a Poisson ideal of $\cA$, which is generated
  (as an ideal) by a single element. If the quotient $\cB:=\cA/I$ is an integral domain then the PBW theorem holds
  for $\cB$.
\end{thm}

As we said above, an important ingredient in our proof is a new construction of the Poisson enveloping algebra
$U(\cB)$ of $\cB:=\cA/I$ in terms of the Poisson enveloping algebra $U(\cA)$ of $\cA$.  In the case of the above
singular example, the quadratic polynomial $X_1^2+X_2^2+X_3^2$ generates a Poisson ideal $I$ of
$\cA:=\bbC[X_1,X_2,X_3]$, equipped with the Poisson bracket (\ref{eq:pb_intro}), so the theorem applies. Here,
$$
  \Omega(\cB)=\frac{\cB\diff X_1+\cB\diff X_2+\cB\diff X_3}{\langle X_1\diff X_1+X_2\diff X_2+X_3\diff X_3\rangle}
$$
and the theorem says that the graded Poisson enveloping algebra of $\cB$ is isomorphic to the symmetric algebra
$\Sym_\cB\Omega(\cB)$. For the precise definition of the isomorphism, which is given by the PBW map, see Section
\ref{par:PBW-map}.

The PBW theorem has important applications to deformation theory and to Poisson and Hochschild (co-) homology; we
will discuss this in a future publication.

The structure of the paper is the following. After quickly recalling the definition of a Poisson algebra and of a
Poisson module over a Poisson algebra we will give the definition of a Poisson enveloping algebra and show that
modules over the latter algebra are in one-to-one correspondence with Poisson modules over the underlying Poisson
algebra. We discuss a few examples of Poisson enveloping algebras of increasing complexity: the ones corresponding
to a null Poisson bracket, to a polynomial algebra, to a general Poisson algebra and to a quotient of a general
Poisson algebra. The latter case is important for Section 3, in which we discuss the PBW theorem. First we
introduce the PBW map and state what it means for a Poisson algebra to satisfy the PBW theorem. We pick up our
list of examples again and show at the end of the section our main result (Theorem \ref{thm:intro_main}). We finish
the paper with some examples and consequences.


\smallskip

In this paper, all rings are assumed to be unitary and all ring morphisms are assumed to preserve the
unit. Similarly, without the adjectives \emph{Lie} or \emph{Poisson}, the word \emph{algebra} stands for an
associative algebra with unit and every {algebra morphism} preserves the unit. Let $R$ be a commutative ring. For
an $R$-module $M$ we denote the tensor algebra of $M$ by $T_R(M)$ or $T(M)$ and the symmetric algebra by
$\Sym_R(M)$ or $\Sym(M)$; both are graded $R$-algebras, with the latter being commutative. For a Lie algebra $\fg$
over $R$, its universal enveloping algebra is denoted by $U_{Lie}(\fg)$. For any algebra $U$ over~$R$ we denote by
$U_L$ the corresponding Lie algebra over $R$, where the bracket is defined by the commutator in $U$: $[u,v]:=uv-vu$
for all $u,v\in U=U_L$. For a graded (resp.\ filtered) algebra $U$ we will denote the factor of $U$ consisting of
all homogeneous elements of degree $n$ by $U^n$ (resp.\ the submodule of all elements of filtered degree at most
$n$ by $U_n$). Without any further specification all algebras are $R$-algebras and $\otimes$ stands for $\otimes_R$.

%

\section{Poisson enveloping algebras}
In this section, we first recall the definition of a Poisson module and of a Poisson enveloping algebra. We give a
few constructions of the Poisson enveloping algebra for the cases of a Lie-Poisson or, more generally, a polynomial
Poisson algebra, the main construction being the construction of the Poisson enveloping algebra of a quotient of a
Poisson algebra by a Poisson ideal. These constructions will turn out to be very useful in the next section, when
we study the PBW theorem for Poisson algebras.

\subsection{The enveloping algebra of a Poisson algebra}
Recall that a \emph{Poisson algebra} (over $R$) is an $R$-module $\cA$ equipped with two multiplications
$(F,G)\mapsto F\cdot G$ and $(F,G)\mapsto\pb{F,G}$, such that
\begin{enumerate}
  \item [(1)] $(\cA,\cdot)$ is a commutative algebra (over $R$);
  \item [(2)] $(\cA,\PB)$ is a Lie algebra (over $R$);
  \item [(3)] The two multiplications are compatible in the sense that the following derivation property is
    satisfied:
  \begin{equation}\label{eq:Leibniz}
     \pb{a_1\cdot a_2,a_3}=a_1\cdot \pb{a_2,a_3}+a_2\cdot \pb {a_1,a_3}\;,
  \end{equation}
  where $a_1,a_2$ and $a_3$ are arbitrary elements of $\cA$.
\end{enumerate}
The bilinear map $\PB$ is called the \emph{Poisson bracket} (of $\cA$). When dealing with the product in a Poisson
algebra, we will always write $a_1a_2$ for $a_1\cdot a_2$. \emph{Morphisms} of Poisson algebras are linear maps
which are both morphisms of algebras and of Lie algebras. A \emph{Poisson ideal} $I$ of $\cA$ is a submodule which
is both an ideal and a Lie ideal of~$\cA$; the quotient $\cB:=\cA/I$ then has a unique Poisson structure, making
the canonical surjection $\pi:\cA\to\cB$ into a morphism of Poisson algebras.

The following three examples will be discussed several times in what follows.
\begin{example}\label{ex:symplectic_poisson}
Let $M$ be an $R$-module and let $\sigma$ be a skew-symmetric bilinear form on $M$. On $\Sym(M)$ a Poisson bracket
is defined setting $\pb{x,y}:=\sigma(x,y)$, for all $x,y\in M$, and extending $\PB$ to a biderivation of
$\Sym(M)$. Explicitly, this yields for monomials $\ux=x_1x_2\dots x_k$ and $\uy=y_1y_2\dots y_\ell$ of $\Sym(M)$,
\begin{equation}\label{eq:symplectic_poisson}
  \pb{\ux,\uy}=\sum_{i=1}^k\sum_{j=1}^\ell\omiet xi\,\omiet yj\,\sigma(x_i,y_j)\;,
\end{equation}%
where $\omiet x i$ stands for the monomial $\ux$ with $x_i$ omitted. When $M=V$ is a finite-dimensional real vector
space and $\sigma$ is a non-degenerate skew-symmetric bilinear form on $V$, then $(V,\sigma)$ is a symplectic
vector space and the above Poisson bracket yields on the algebra $C^\infty(V)$ a Poisson bracket, which is
precisely Poisson's original bracket (see \cite[Ch.~6]{PLV}).
\end{example}
\begin{example}\label{ex:affine}
Let $(\fg,\LB)$ be a Lie algebra over $R$ and let $\sigma$ be a 2-cocycle in the trivial Lie algebra cohomology of
$(\fg,\LB)$. The latter means that $\sigma$ is a skew-symmetric bilinear form on $\fg$, such that
$$
  \sigma([x,y],z)+\sigma([y,z],x)+\sigma([z,x],y)=0\;,
$$
for all $x,y,z\in\fg$.  A Poisson bracket is defined on $\Sym(\fg)$ by setting
$\pb{x,y}_\sigma:=\lb{x,y}+\sigma(x,y)$ for all $x,y\in\fg$, and extending $\PB_\sigma$ to a biderivation of
$\Sym(\fg)$. Explicitly, it is given as in (\ref{eq:symplectic_poisson}), with $\sigma(x_i,y_j)$ replaced by
$\lb{x_i,y_j}+\sigma(x_i,y_j)$. When the Lie bracket $\LB$ is the trivial bracket, the present example reduces to
Example~\ref{ex:symplectic_poisson}. When $\sigma$ is trivial, $\PB_\sigma$ is a linear Poisson structure, usually
referred to as a \emph{Lie-Poisson structure} and $(\Sym(\fg),\PB)$ is called a \emph{Lie-Poisson algebra}; in
general, we refer to it as a \emph{modified Lie-Poisson algebra} (see \cite[Ch.\ 7]{PLV}).
\end{example}
\begin{example}\label{exa:Nambu}
Let $P$ and $Q$ be two polynomials in three variables. They define a Poisson structure on $R[X_1,X_2,X_3]$ by
setting
\begin{equation*}
  \pb{X_1,X_2}:=Q\pp P{X_3}\;,\quad   \pb{X_2,X_3}:=Q\pp P{X_1}\;,\quad   \pb{X_3,X_1}:=Q\pp P{X_2}\;,
\end{equation*}%
which is again extended to a biderivation of $R[X_1,X_2,X_3]$. Notice that $P$ is a \emph{Casimir} of this Poisson
structure, i.e., it belongs to the center of $\PB$. The above Poisson structure on $R[X_1,X_2,X_3]$ is called a
\emph{Nambu-Poisson structure} (see \cite[Ch.\ 8.3]{PLV}).
\end{example}
\begin{example}\label{ex:poisson_poly}
Suppose that $M$ is an $R$-module and that its symmetric algebra $\Sym(M)$ is equipped with a skew-symmetric
biderivation $\PB$, satisfying the Jacobi identity for all triplets of elements from~$M$. Then~$\PB$ satisfies the
Jacobi identity for all triplets of elements from $\Sym(M)$, hence makes $\Sym(M)$ into a Poisson algebra. Such and
algebra is called a \emph{polynomial Poisson algebra} (see \cite[Ch.\ 1.4, 8.1]{PLV}).
\end{example}
Let $\cA$ be a Poisson algebra (over $R$). A \emph{Poisson module} over $\cA$ is an $R$-module~$E$ which is both a
module and a Lie module over $\cA$, satisfying supplementary derivation (Leibniz) rules (see
\cite{cherkashin,MR1683858}). To be precise, $\cA$ is equipped with two maps $\a_E,\b_E:\cA\to\End(E)$, such that,
for all $a_1,a_2\in\cA$,
\begin{enumerate}
  \item [(1)] $\a_E(a_1 a_2)=\a_E(a_1)\a_E(a_2)$ 
  \item [(2)] $\b_E(\pb{a_1, a_2})=\b_E(a_1)\b_E(a_2)-\b_E(a_2)\b_E(a_1)$
  \item [(3)] $\a_E(\pb{a_1,a_2})=\a_E(a_1)\b_E(a_2)-\b_E(a_2)\a_E(a_1)$,
  \item [(4)] $\b_E(a_1 a_2)=\a_E(a_1)\b_E(a_2)+\a_E(a_2)\b_E(a_1)$.
\end{enumerate}
In the right hand side of these formulas, the product is composition of elements of $\End(E)$. Item (1) resp.\ (2)
says that $\a_E$ is an algebra morphism, resp.\ that $\b_E$ is a Lie algebra morphism; item (4) says that $\b_E$ is
an $\a_E$-derivation.

Examples of Poisson modules include $\cA$ itself and any of its powers, any Poisson ideal of $\cA$, the dual of
$\cA$ and so on. With the natural notion of morphism between Poisson modules (where one asks that the morphism is
both a morphism of modules and of Lie modules), the Poisson modules over $\cA$ form a category which, as we will
see later, is an Abelian category.

We are now ready for defining the notion of a Poisson enveloping algebra.
\begin{defn}\label{def:PEA}
  Let $(\cA,\cdot,\PB)$ be a Poisson algebra (over $R$). A \emph{Poisson enveloping algebra} for $\cA$ is an
  algebra $U$, equipped with two maps:
\begin{enumerate}
  \item[(1)] An algebra morphism $\a:(\cA,\cdot)\to U$,
  \item[(2)] A Lie algebra morphism $\b:(\cA,\PB)\to U_L$, 
\end{enumerate}
such that, for any $a_1,a_2\in\cA$,
\begin{enumerate}
  \item[(3)] $\a(\pb{a_1,a_2})=\a(a_1).\b(a_2)-\b(a_2).\a(a_1),$
  \item[(4)] $\b(a_1a_2)=\a(a_1).\b(a_2)+\a(a_2).\b(a_1),$
\end{enumerate}
and such that the following universal property holds: if $U'$ is any algebra and $\a':(\cA,\cdot)\to U'$ and
$\b':(\cA,\PB)\to U'_L$ are any algebra (resp.\ Lie algebra) morphisms, satisfying the following properties: for any
$a_1,a_2\in\cA$,
\begin{enumerate}
  \item[(3')] $\a'(\pb{a_1,a_2})=\a'(a_1).\b'(a_2)-\b'(a_2).\a'(a_1),$
  \item[(4')] $\b'(a_1a_2)=\a'(a_1).\b'(a_2)+\a'(a_2).\b'(a_1),$
\end{enumerate}
then there exists a unique algebra morphism $\gamma:U\to U'$, such that
$$
  \gamma\circ \a=\a',\qquad \gamma\circ \b=\b'\;.
$$
The two equalities are summarized in the following commutative diagram:
\begin{equation}\label{dia:poisson_envel}
  \begin{tikzpicture}
  \matrix (m) [matrix of math nodes,row sep=3em,column sep=4em,minimum width=2em]
  {
     U \\
     \cA & U' \\};
  \path[-stealth]
    (m-2-1) edge node [left] {$\a,\b$} (m-1-1)
            edge node [below] {$\a',\b'$} (m-2-2)
    (m-1-1) edge [dashed,->] node [right] {$\gamma$} (m-2-2);
  \end{tikzpicture}
\end{equation}
\end{defn}
Notice that $U$ is not only an $R$-algebra, but is also in a natural way an $\cA$-module, as we may define, for
$a\in \cA$ and for $u\in U$, $a\cdot u:=\a(a)u$, which we often write simply as $au$; in general, we often
identify $a\in\cA$ with $\a(a)\in U(\cA)$, which is without danger because the algebra morphism $\a$ is
always an injection (see \cite[p.\ 198]{rinehart}, and also \cite[Prop.\ 2.2]{oh_park_shin}). Notice that the
algebra morphism $\gamma:U\to U'$ in diagram (\ref{dia:poisson_envel}) is a morphism of $\cA$-modules, when $U$ and
$U'$ are viewed as $\cA$-modules.
\begin{thm}
  Let $(\cA,\cdot,\PB)$ be a Poisson algebra (over $R$). There exists a Poisson enveloping algebra for $\cA$ and it
  is unique up to isomorphism: if $(U,\a,\b)$ and $(U',\a',\b')$ are two Poisson enveloping algebras for $\cA$,
  then there exists an algebra isomorphism $\gamma:U\to U'$, such that $\gamma\circ\a=\a'$ and
  $\gamma\circ\b=\b'$. The Poisson enveloping algebra of $\cA$, which is unique up to isomorphism, is denoted by
  $U(\cA)$ and its accompanying maps are denoted by $\a$ and $\b$ (or by $\a_\cA$ and $\b_\cA$ when more than one
  Poisson enveloping algebra is considered).
\end{thm}
Uniqueness of the Poisson enveloping algebra is clear. A few different existence proofs can be found in
\cite{rinehart,huebsch,MR1683858}. We will give in the next subsections a few alternative constructions for the
cases which we will consider in the next section. In the present subsection, we only treat the example of a Poisson
algebra $\cA$ whose Poisson bracket is the zero bracket. This case is quite simple, but very instructive, as it
will provide the natural candidate for the source of the PBW map and give a first instance of a Poisson algebra
which satisfies the PBW theorem (stated and treated in general in Section \ref{sec:PBW}). For any algebra $\cA$ we
denote by $\Omega(\cA)$ the $\cA$-module of K\"ahler differentials on $\cA$ (see \cite[Ch.\ 16]{eisenbud}).
\begin{prop}\label{prop:null}
  Let $\cA$ be any algebra which we make into a Poisson algebra by adding the zero Poisson bracket. Denote by
  $\a:\cA\to \Sym_\cA(\Omega(\cA))$ the canonical inclusion map and let $\b:=\diff:\cA\to
  \Sym_\cA(\Omega(\cA))$. The triplet $(\Sym_\cA(\Omega(\cA)),\a,\b)$ is a Poisson enveloping algebra for $\cA$.
\end{prop}
\begin{proof}
Since the Poisson bracket on $\cA$ is null and since $\Sym_\cA(\Omega(\cA))$ is commutative, the verification of
properties (1) -- (3) in Definition \ref{def:PEA} is immediate; (4) is just the relation $\diff(a_1a_2)=a_1\diff
a_2+a_2\diff a_1$ (which holds in $\Omega(\cA)$), rewritten in terms of $\a$ and $\b$. Suppose now that $U'$ is any
algebra and that $\a',\b':\cA\to U'$ are any algebra (resp.\ Lie algebra) morphisms, satisfying (3') and (4') in
Definition \ref{def:PEA}. This means in particular that all elements in the image of $\a'$ and $\b'$ commute. If
there exists an algebra morphism $\gamma:\Sym_\cA(\Omega(\cA))\to U'$, such that $\gamma\circ \a=\a'$ and
$\gamma\circ \b=\b',$ then it is given by
\begin{eqnarray*}
  \gamma(a\diff a_1\diff a_2\dots\diff a_k)=\a'(a).\b'(a_1).\b'(a_2)\dots \b'(a_k)\;,
\end{eqnarray*}
for $a,a_1,\dots,a_k\in\cA$.  This shows that $\gamma$ is unique, if it exists. Clearly, $\gamma$ is well-defined
by this formula and satisfies $\gamma\circ \a=\a'$ and $\gamma\circ \b=\b'.$ Finally, $\gamma$ is an algebra
morphism because all elements in the image of $\a'$ and $\b'$ commute.
\end{proof}

Let $E$ be a Poisson module over $\cA$, with structure maps $\a'$ and $\b'$. In view of the universal property of
the Poisson enveloping algebra, there exists an algebra morphism $\gamma:U(\cA)\to\End(E)$ as in
(\ref{dia:poisson_envel}) (with $U'=\End(E)$), in particular $E$ has a natural structure of
$U(\cA)$-module. Conversely, composition with $\a$ and $\b$ transforms any $U(\cA)$ module into a Poisson module
over~$\cA$. The upshot of this natural (and functorial) correspondence is that the category of Poisson modules over
$\cA$ is equivalent to the category of modules over $U(\cA)$. It follows that the category of Poisson modules over
a given Poisson algebra is an Abelian category, as we announced earlier.

%
%
%

\subsection{The Poisson enveloping algebra of a polynomial Poisson algebra}\label{par:envel_poly}
In this subsection we give a new construction of the Poisson enveloping algebra of a polynomial Poisson algebra,
which we also specialize to the case of a Lie-Poisson algebra. For doing this, we use smash product algebras, which
are constructed from module algebras, two notions which we first recall (for more details on these notions and for
proofs, see \cite{montgomery}, for example).

\subsubsection{Module algebras and smash product algebras}
Let $H$ be a Hopf algebra and let $\cA$ be an algebra (both over~$R$). One says that $\cA$ is a \emph{(left)
$H$-module algebra} if $\cA$ has the structure of a left $H$-module, with the following properties: for all $u\in
H$ and for all $a_1,a_2\in\cA$,
\begin{enumerate}
  \item[(1)] $u\cdot(a_1a_2)=\sum_{(u)}(u_{(1)}\cdot a_1)(u_{(2)}\cdot a_2)$;
  \item[(2)] $u\cdot 1=\epsilon(u) 1$.
\end{enumerate}
In (1) we have used Sweedler's notation, i.e., we have written the coproduct of $u\in H$ as $\Delta(u)=\sum_{(u)}
u_{(1)}\otimes u_{(2)}$. Also, $\epsilon$ denotes the counit of $H$. The \emph{smash product algebra} of $\cA$ by
$H$, denoted $\cA\# H$ is as an $R$-module $\cA\otimes H$, with elements denoted by $a\# u$, and with product
defined for all $a_1,a_2\in \cA$ and $u,v\in H$ by
\begin{equation}\label{eq:def_smash}
  (a_1\# u)\odot(a_2\# v) :=\sum_{(u)} a_1(u_{(1)}\cdot a_2)\# u_{(2)}v\;.  
\end{equation}%
This product is associative with unit $1\#1$. For $a_1,a_2\in \cA$ and $u,v\in H$ it follows from definition
(\ref{eq:def_smash}) that $(a_1\#1)\odot(a_2\#1)=a_1(1\cdot a_2)\#1=a_1a_2\#1,$ so that the inclusion map
$\imath_\cA:\cA \rightarrow \cA\# H$ is a morphism of algebras. It can be used to define an $\cA$-module structure
on $\cA\# H$ by setting, for $a_1,a_2\in\cA$ and $u\in H$,
\begin{equation*}
  a_1\cdot(a_2\#u):=\imath_\cA(a_1)\odot(a_2\#u)=a_1a_2\#u\;.
\end{equation*}%
In the sequel, we write $a_1(a_2\#u)$ for $a_1\cdot(a_2\#u)$. It also follows from the definitions that
\begin{eqnarray}\label{eq:hopf_sub}
  (a\#u)\odot(1\#v)&=&\sum_{(u)}a(u_{(1)}\cdot 1)\#u_{(2)}v=\sum_{(u)}a\epsilon(u_{(1)})\#u_{(2)}v\nonumber\\
  &=&a\#\left(\sum_{(u)}\epsilon(u_{(1)})u_{(2)}\right)v=a\#uv\;,
\end{eqnarray}
where we used in the last equality that $\epsilon$ is the counit of $H$. This shows in particular that the
inclusion map $\imath_H:H \rightarrow \cA\# H$ is also a morphism of algebras. Notice that every element of $\cA\#
H$ of the form $a\#u$ can be written as the product of an element of $\Im(\imath_\cA)$ with an element of
$\Im(\imath_H)$, namely for all $a\in\cA$ and $u\in H$,
\begin{equation}\label{for:smash_generated}
  (a\#1)\odot(1\#u)=a\#u\;.
\end{equation}%
It leads, in view of the above properties, to a simple proof that $1\#1$ is the unit of $\cA\# H$, as we said
above. We give two typical examples; they will be used later, besides others which will be introduced as we need
them.
\begin{example}\label{exa:smash_poly}
Let $M$ be an $R$-module and let $\PB$ be a Poisson bracket on $\Sym(M)$, making it into a polynomial Poisson
algebra. It is well-known (see \cite[Ch.\ 3]{kassel_book}) that the tensor algebra $T(M)$ has a natural structure
of a Hopf algebra, where the comultiplication $\Delta:T(M)\to T(M)\otimes T(M)$ is the unique algebra morphism,
given for $x\in M$ by $\Delta(x)= 1\otimes x+x\otimes 1$ and the counit $\epsilon:T(M)\to R$ picks the constant
(degree zero) term of a tensor. Using the Poisson bracket, $\Sym(M)$ becomes a $T(M)$-module algebra upon setting,
for $x_1\otimes x_2\otimes\cdots\otimes x_k\in T(M)$ and $a\in \Sym(M)$,
\begin{equation}\label{eq:mod_alg_tens}
  (x_1\otimes x_2\otimes\cdots\otimes x_k) \cdot a:=\pb{x_1,\pb{x_2,\dots\pb{x_k,a}\dots}}\;.
\end{equation}%
It is understood that this definition specializes for $k=0$ to $1\cdot a:=a$. Since both $\Delta$ and $\epsilon$
are algebra morphisms, and since items (1) and (2) in the above definition of a module algebra are obviously
satisfied for $u=1$, it suffices to check them for $u=x\in M$ (and $a_1,a_2\in \Sym(M)$). Since $\PB$ is a
derivation in each argument, we have
\begin{equation*}
  x\cdot(a_1a_2)=\pb{x,a_1a_2}=a_1\pb{x,a_2}+\pb{x,a_1}a_2=(1\cdot a_1)(x\cdot a_2)+(x\cdot a_1)(1\cdot a_2)\;,
\end{equation*}%
which proves (1), since $\Delta(m)=1\otimes m+m\otimes 1$. Also, (2) holds because $m\cdot 1=\pb{m,1}=0$ and
$\epsilon(m)=0$. This shows that $\Sym(M)$ is a $T(M)$-module algebra. We can therefore form the smash product
algebra $\Sym(M)\# T(M)$. According to definition (\ref{eq:def_smash}), the product in $\Sym(M)\# T(M)$ is given,
for $a_1,a_2\in \Sym(M)$ and $x\in M$ and $u\in T(M)$ by
\begin{equation}\label{eq:def_smash_tns}
  (a_1\# x)\odot(a_2\# u) =a_1a_2\#(x\otimes u)+a_1\pb{x,a_2}\#u\;.  
\end{equation}%

\end{example}

\begin{example}\label{ex:smash_general}
Let $\cA$ be an arbitrary Poisson algebra. We show that $\cA$ is a $U_{Lie}(\cA)$ module algebra. To do this, we
first recall that the standard Hopf algebra structure of $U_{Lie}(\cA)$ is induced by the Hopf algebra structure on
$T(\cA)$, recalled in the previous example (see \cite[Ch.\ 5]{kassel_book}). It means that $\Delta: U_{Lie}(\cA)\to
U_{Lie}(\cA)\otimes U_{Lie}(\cA)$ is the unique algebra morphism which is defined for all $a\in\cA\subset
U_{Lie}(\cA)$ by $\Delta(a):=1\otimes a+a\otimes 1$. Consider the linear map $\X:\cA\to\End(\cA)$, defined by
$a\mapsto \X_a=\pB{a}$. In view of the Jacobi identity, it makes $\cA$ into a Lie module over $\cA$, hence makes
$\cA$ into a module over $U_{Lie}(\cA)$. To check that it makes $\cA$ into a (left) $U_{Lie}(\cA)$ module algebra,
it suffices to verify (1) and (2) in the above definition of a module algebra for $u\in\cA$ and for $u\in R$ (that
is for $u\in U_{Lie}(\cA)$ of degree at most 1). For example, when $u\in\cA$, so that $\Delta(u)=1\otimes
u+u\otimes 1$, then $u\cdot (a_1a_2)=\pb{x,a_1a_2}$, while
\begin{equation*}
  \sum_{(u)}(u_{(1)}\cdot a_1)(u_{(2)}\cdot a_2)=(1\cdot a_1)(u\cdot a_2)+(u\cdot a_1)(1\cdot a_2)
   =a_1\pb{u,a_2}+a_2\pb{u,a_1}\;,
\end{equation*}%
which is the same as $\pb{u,a_1b_1}$ because $\PB$ is a biderivation. The other verifications are even
simpler. Thus, $\cA$ is a left $U_{Lie}(\cA)$-module algebra and we can form the smash product $\cA\# U_{Lie}(\cA)$
of $\cA$ by $U_{Lie}(\cA)$. For future use, let us point out that the product in $\cA\# U_{Lie}(\cA)$ is given, for
$a_1,a_2,a_3\in\cA$ and $u\in U_{Lie}(\cA)$ by
\begin{equation}\label{eq:def_smash_2}
  (a_1\# a_3)\odot(a_2\# u) =a_1a_2\#a_3.u+a_1\pb{a_3,a_2}\#u\;,
\end{equation}%
where $a_3.u$ stands for the product of $a_3$ and $u$ in $U_{Lie}(\cA)$.
\end{example}
%

%
%
\subsubsection{The Poisson enveloping algebra of a (modified) Lie-Poisson algebra as a smash product algebra}
\label{par:envel_lie-poisson}
We show in this paragraph that the Poisson enveloping algebra of the modified Lie-Poisson algebra
$\Sym_\sigma(\fg)=(\Sym(\fg),\PB_\sigma)$ (where $\fg$ is a Lie algebra and $\sigma$ is a 2-cocycle in the trivial
Lie algebra cohomology of~$\fg$, see Example \ref{ex:affine}) is the smash product algebra
$\Sym(\fg)\#U_{Lie}(\fg)$, with accompanying maps $\a$ and $\b$ which will be defined below.

First, we need to explain how we turn $\Sym(\fg)$ into a $U_{Lie}(\fg)$-module algebra. The construction is very
similar to the one given in Example \ref{ex:smash_general}: $\Sym(\fg)$ is a Lie module over $\fg$, when setting
$x\cdot a:=\pb{x,a}_\sigma$, for $x\in\fg$ and $a\in \Sym(\fg)$, so that $\Sym(\fg)$ is a module over
$U_{Lie}(\fg)$ and one verifies like in Example \ref{ex:smash_general} that $\Sym(\fg)$ is a $U_{Lie}(\fg)$-module
algebra. We can therefore form the smash product algebra $\Sym(\fg)\#U_{Lie}(\fg)$. By construction, the product in
this algebra is given, for $a_1,a_2\in \Sym_{\sigma}(\fg)$ and $x\in\fg\subset U_{Lie}(\fg)$ and $u\in
U_{Lie}(\fg)$ by
\begin{equation}\label{for:smash_lie-poisson}
  (a_1\# x)\odot(a_2\# u) =a_1a_2\#x.u+a_1\pb{x,a_2}_\sigma\#u\;.  
\end{equation}%
%
%
The maps $\a$ and $\b$ are defined by
\begin{equation*}
  \begin{array}{lcccl}
    \a&:&\Sym_\sigma(\fg)&\to&\Sym_\sigma(\fg)\# U_{Lie}(\fg)\\
    & & a&\mapsto&a\#1\;,
  \end{array}
\end{equation*}
\begin{equation*}
  \begin{array}{lcccl}
    \b&:&\Sym_\sigma(\fg)&\to&\Sym_\sigma(\fg)\# U_{Lie}(\fg)\\
    & & \ux&\mapsto&\sum_i\omiet xi\# x_i\;.
  \end{array}
\end{equation*}
Notice that $\b$ can be defined as the unique $\a$-derivation such that $\b(x)=1\#x$ for all $x\in\fg$.
\begin{prop}\label{thm:envel_lie-poisson}
  $(\Sym(\fg)\#U_{Lie}(\fg),\a,\b)$ is a Poisson enveloping algebra of the Lie-Poisson algebra
  $\Sym_\sigma(\fg)$.
\end{prop}
\begin{proof}
We first verify items (1) -- (4) in Definition \ref{def:PEA}. We know that $\a$ is an algebra morphism and that
$\b$ is an $\a$-derivation, which is the content of~(1) and (4). We move to item (3): for a monomial
$\ux\in \Sym_\s(\fg)$ and for any $a\in \Sym_\s(\fg)$ we have
\begin{eqnarray*}
  \lb{\a(a),\b(\ux)}
  &=&a\#1\odot\sum_i\omiet xi\#x_i-\sum_i\omiet xi\#x_i\odot a\#1\\
  &=&\sum_ia\,\omiet xi\#x_i-\sum_i\omiet xi a\#x_i-\sum_i\omiet xi\pb{x_i,a}_\sigma\#1\\
  &=&\pb{a,\ux}_\sigma\#1=\a(\pb{a,\ux}_\sigma)\;.
\end{eqnarray*}
%
%
We still need to prove (2). In view of the items just proved, it is easily shown by recursion that it suffices to
prove that $\b(\pb{x,y}_\sigma)=\lb{\b(x),\b(y)}$ for all $x,y\in\fg$. Since
$\pb{x,y}_\sigma=\lb{x,y}+\sigma(x,y)$ and since $\b$ is null on constants, we have
\begin{eqnarray*}
  \b(\pb{x,y}_\sigma)&=&\b(\lb{x,y})=1\#\lb{x,y}=1\#(x.y-y.x)\\
      &=&(1\#x)\odot(1\#y)-(1\#y)\odot(1\#x)=\lb{\b(x),\b(y)}\;.
\end{eqnarray*}
We now prove the universal property of $\Sym(\fg)\#U_{Lie}(\fg)$. Let $U'$ be any algebra and suppose that we are
given any algebra morphism $\a':\Sym_\sigma(\fg)\to U'$ and any Lie algebra morphism $\b':\Sym_\sigma(\fg)\to
U'_L$, satisfying properties (3') and (4') (of Definition \ref{def:PEA}). We show that there is a unique algebra
morphism $\gamma:\Sym(\fg)\# U_{Lie}(\fg)\to U'$ such that $\gamma\circ\a=\a'$ and $\gamma\circ\b=\b'$. Since every
element $a\#u$ of $\Sym(\fg)\# U_{Lie}(\fg)$ can be written as the product of an element of $\Im(\a)$ with a
product of elements of $\Im(\b)$, the morphism $\gamma$ is unique, if it exists; moreover, it leads to the formulas
\begin{eqnarray*}
  \gamma(a\#1)&=&\a'(a)\;,\\
  \gamma(a\#(x_1.x_2\dots x_k))&=&\a'(a).\b'(x_1).\b'(x_2)\dots \b'(x_k)\;,
\end{eqnarray*}
where $a\in\Sym_\s(\fg)$ and $x_1,\dots,x_k\in\fg$. In order to show that the map $\gamma$ is well-defined by this
formula, we use the $R$-linear map $\gamma':\Sym_\sigma(\fg)\#T(\fg)\to U'$, defined by
\begin{equation*}
  \gamma'(a\otimes x_1\otimes\cdots\otimes x_k):=\a'(a).\b'(x_1)\dots\b'(x_k)\;,  
\end{equation*}%
where $a\in\Sym_\sigma(\fg)$ and $x_1,\dots,x_k\in\fg$. For $X=x_1\otimes x_2\otimes\cdots\otimes x_k$ and
$Y=y_1\otimes y_2\otimes\cdots \otimes y_\ell$ in~$T(\fg)$, and for $x,y\in\fg$, we have
\begin{eqnarray*}
  \lefteqn{{\gamma}'(a\otimes X\otimes (x\otimes y-y\otimes x-[x,y])\otimes Y)=}\hskip 12cm\\
  \a'(a).\b'(x_1)\dots\beta'(x_k).([\b'(x),\b'(y)]-\b'([x,y])).\b'(y_1)\dots\beta'(y_\ell)=0\;,
\end{eqnarray*}
because $\beta'$ is a Lie morphism. It follows that the $R$-linear map $\gamma: \Sym_\sigma(\mathfrak{g}) \otimes
U_{Lie}(\mathfrak{g})\to U'$ is well-defined.  We need to show that $\gamma$ is an algebra morphism, i.e.\ that
$\gamma((a_1\#u)\odot(a_2\#v))= \gamma(a_1\#u).\gamma(a_2\#v)$ for all $a_1,a_2\in \Sym(\fg)$ and for all
homogeneous elements $u,v$ of $U_{Lie}(\fg)$. We do this by recursion on the filtered degree $k$ of $u$ and we
write $v=v_1.v_2\dots v_\ell$, where all $v_i$ belong to $\fg$. For $u=1$, we have
\begin{equation*}
  \gamma((a_1\#1)\odot(a_2\#v))=\a'(a_1a_2).\b'(v_1)\dots\b'(v_\ell)=\gamma(a_1\#1).\gamma(a_2\#v)\;.
\end{equation*}%
We next take $u=x\in\fg\subset U_{Lie}(\fg)$. Then, using (3'),
\begin{eqnarray*}
  \gamma((a_1\#x)\odot(a_2\#v))
  &=&\gamma(a_1a_2\#x. v+a_1\pb{x,a_2}\#v)\\
  &=&\a'(a_1).\left(\a'(a_2).\b'(x)-\a'(\pb{a_2,x})\right).\b'(v_1)\dots\b'(v_\ell)\\
  &=&\a'(a_1).\b'(x).\a'(a_2).\b'(v_1)\dots\b'(v_\ell)\\
  &=&\gamma(a_1\#x)\,.\gamma(a_2\#v)\;.
\end{eqnarray*}
Suppose now that $\gamma((a_1\#u)\odot(a_2\#v))=\gamma(a_1\#u).\gamma(a_2\#v)$ holds for any $u$ of degree at most
$k$. Then, using the associativity of $\odot$, the recursion hypothesis and (\ref{eq:hopf_sub}),
\begin{eqnarray*}
{\gamma((a_1\#u.x)\odot(a_2\#v))}
  &=&\gamma(a_1\#u\odot(1\#x\odot a_2\#v))\\
  &=&\gamma(a_1\#u)\,.\gamma(1\#x\odot a_2\#v)\\
  &=&\gamma(a_1\#u)\,.\gamma(1\#x)\,.\gamma(a_2\#v)\\
  &=&\gamma(a_1\#u\odot 1\#x)\,.\gamma(a_2\#v)\\
  &=&\gamma(a_1\#u.x)\,.\gamma(a_2\#v)\;.
\end{eqnarray*}
This shows that the formula also holds for $u$ of degree at most $k+1$, and hence that $\gamma$ is an algebra
morphism. Finally, we need to check that $\gamma \circ\a=\a'$ and $\gamma \circ\b=\b'$. The first equality is
immediate from the above explicit formula for $\gamma $, so we only prove the second one. For any monomial $\ux\in
\Sym_\sigma(\fg)$,
\begin{eqnarray*}
  \gamma (\b(\ux))=\sum_i \gamma (\omiet xi\#x_i)=\sum_i\a'(\omiet xi).\b'(x_i)=\b'(\ux)\;.
\end{eqnarray*}
The last equality is valid because $\b'$ is an $\a'$-derivation.
\end{proof}

\subsubsection{The Poisson enveloping algebra of a (modified) Lie-Poisson algebra as a (modified) Lie enveloping
  algebra}\label{par:envel_lie-poisson_2}
We give in this paragraph a different description of the Poisson enveloping algebra of a (modified) Lie-Poisson
algebra. In the unmodified case, the result is that for any Lie algebra $\fg$, the Lie enveloping algebra of a
certain double $\fg^+$ of $\fg$ (known as a Takiff algebra, see \cite{Takiff}) is a Poisson enveloping algebra of
the Lie-Poisson algebra $\Sym(\fg)$. In the modified case, the same result holds, upon using the notion of a
modified Lie enveloping algebra, also known as a Sridharan algebra (\cite{Sridharan}).

Suppose, as in the previous paragraph, that $\fg$ is a Lie algebra and that $\sigma$ is a 2-cocycle in the trivial
Lie algebra cohomology of~$\fg$. Let us denote by $\fg^0$ the abelian Lie algebra, whose underlying module is
$\fg$.  Consider $\fg^+:=\fg^0\oplus\fg$, in which $\fg^0$ and $\fg$ are naturally embedded. For $x\in\fg$ we will
write $x^0$, respectively $x^1$, for its canonical image in $\fg^0$, respectively in~$\fg$, viewed as a subspace of
$\fg^+$. Thus, we can write every element $x^+$ of $\fg^+$ uniquely as $x^+=y^0+z^1$, with $y,z\in\fg$.  A Lie
bracket is defined on $\fg^+$ by
$$
  [y^0,z^0]^+=0,\quad [y^0,z^1]^+=[y,z]^0,\quad [y^1,z^1]^+=[y,z]^1\;,
$$
where $y,z\in\fg$. The Lie algebra $\fg^+$ is a semi-direct product of $\fg^0$ and $\fg$: for elements
$x_1^+=y_1^0+z_1^1$ and $x_2^+=y_2^0+z_2^1$ of $\fg^+$, we have 
\begin{equation*}
  [x_1^+,x_2^+]^+=[y_1,z_2]^0+[z_1,y_2]^0+[z_1,z_2]^1\;.
\end{equation*}%
The cocycle $\s$ becomes a cocycle $\s^+$ of $\fg^+$ upon setting for all $x,y\in\fg$:
$$
  \sigma^+(x^0,y^0):=\sigma^+(x^1,y^1):=0,\ \sigma^+(x^0,y^1):=\sigma^+(x^1,y^0):=\sigma(x,y)\;,
$$
  and extending these definitions by bilinearity. Since
$$
  \sigma^+([x_1^+,x_2^+],x_3^+)=\s([y_1,z_2],z_3)+\s([z_1,y_2],z_3)+\s([z_1,z_2],y_3)
$$  
and since $\s$ is a cocycle, $\s^+$ is indeed a cocycle. The modified Lie enveloping algebra (or Sridharan algebra)
of $\Sym_\s(\fg)$ is given by $U_{Lie,\sigma^+}(\fg^+):=T(\fg)/I_\sigma$ where $I_\sigma$ is the two-sided ideal of
$T(\mathfrak{g})$ generated by all elements of the form
\begin{equation}\label{eq:Sri_ideal}
  \{x^+\otimes y^+-y^+\otimes x^+-[x^+,y^+]^+-\sigma^+(x^+,y^+)\cdot 1\}\;,
\end{equation}
where $x^+,y^+\in\fg^+$.  Let $\iota$ denote the canonical inclusion $\iota: \fg^+\hookrightarrow
U_{Lie,\sigma^+}(\fg^+)$.  For $x\in\fg$, let $\a'(x):=\iota(x^0)$ and $\b'(x):=\iota(x^1)$. For $x,y\in\fg$ we
have, in view of (\ref{eq:Sri_ideal}), that $\iota(x^0)\iota(y^0)-\iota(y^0)\iota(x^0)= \iota([x^0,y^0]^+)
+\sigma^+(x^0\otimes y^0)\cdot 1=0$. We can therefore uniquely extend $\a'$ to an algebra morphism
$$
  \a':\Sym_\sigma(\fg) \to U_{Lie,\sigma^+}(\fg^+)\;.
$$
By a slight abuse of notation, we will also write $\iota(a^0)$ for $\a'(a)$, where $a\in \Sym_\sigma(\fg)$.  As for
$\b'$, it extends uniquely to an $\a'$-derivation
$$
  \b':\Sym_\sigma(\fg) \to U_{Lie,\sigma^+}(\fg^+)\;.
$$
Explicitly, $\b'$ is given for a monomial $\ux\in \Sym_\sigma(\fg)$ by $\b'(\ux)=\sum_j\iota(\omiet x
j)^0\iota(x_j^1)$.
\begin{prop}\label{prp:envel_lie-poisson}
  The modified Lie enveloping algebra $(U_{Lie,\s^+}(\fg^+), \a',\b')$ is a Poisson enveloping algebra of the
  modified Lie-Poisson algebra $\Sym_{\sigma}(\fg)$.
\end{prop} 
\begin{proof}
By construction, $\a'$ is an algebra morphism and $\a'$ and $\b'$ satisfy property (4') of Definition
\ref{def:PEA}. We show that they also satisfy property (3') of the latter definition and that $\b'$ is a morphism
of Lie algebras. For a monomial $\ux\in \Sym_\s(\fg)$ and for any $a\in \Sym_\s(\fg)$ we have
\begin{eqnarray*}
  \lb{\a'(a),\b'(\ux)} &=&\iota(a^0).\sum_j\iota((\omiet x j)^0).\iota(x_j^1)-
                          \sum_j\iota((\omiet x j)^0).\iota(x_j^1).\iota(a^0)  \\
  &=&-\sum_j\iota((\omiet x j)^0).\iota(\pb{x_j,a}_\s^0) =\iota(\pb{a,\ux}_\sigma^0)=\a'(\pb{a,\ux}_\sigma)\;.
\end{eqnarray*}
As in the proof of Proposition \ref{thm:envel_lie-poisson}, $\b'$ is a Lie morphism as soon as it has the Lie
morphism property when applied to elements of $\fg$. Therefore, let $x,y\in\fg$. On the one hand,
$\b'(\pb{x,y}_\sigma)=\b'([x,y]+\sigma(x,y))=\b'([x,y])=\iota([x,y]^1)$, while on the other hand,
$[\b'(x),\b'(y)]=\iota(x^1).\iota(y^1)-\iota(y^1).\iota(x^1)=\iota([x,y]^1)+\sigma^+(x^1,y^1)\cdot
1=\iota([x,y]^1)$. This shows that $\beta'$ is a Lie algebra morphism.

We can now apply the universal property of the Poisson enveloping algebra $(\Sym(\fg)\#U_{Lie}(\fg),\a,\b)$ (see
Proposition \ref{thm:envel_lie-poisson}): there exists a (unique) algebra morphism $\gamma$, making the following
diagram commutative:
\begin{equation*}
  \begin{tikzpicture}
  \matrix (m) [matrix of math nodes,row sep=3em,column sep=4em,minimum width=2em]
  {
     \Sym(\fg)\#U_{Lie}(\fg) \\
     \Sym_\s\fg &U_{Lie,\s^+}(\fg^+)  \\};
  \path[-stealth]
    (m-2-1) edge node [left] {$\a,\b$} (m-1-1)
            edge node [below] {$\a',\b'$} (m-2-2)
    (m-1-1) edge [dashed,->] node [right] {$\gamma$} (m-2-2);
  \end{tikzpicture}
\end{equation*}
In order to show that $\gamma$ is an isomorphism, we construct the inverse map. We use for it the universal
property of the modified Lie enveloping algebra $U_{Lie,\s^+}(\fg^+)$. Denote by $\jmath:\fg^+\to
\Sym(\fg)\#U_{Lie}(\fg)$ the linear map, defined by $\jmath(x^+)=\jmath(y^0+z^1):=\a(y)+\b(z)$. For
$x_1^+=y_1^0+z_1^1$ and $x_2^+=y_2^0+z_2^1$ in $\fg^+$ we have the following three equalities:
\begin{eqnarray*}
  \lb{\jmath(x_1^+),\jmath(x_2^+)}
    &=&\lb{\a(y_1)+\b(z_1),\a(y_2)+\b(z_2)}\\
    &=& \lb{\a(y_1),\b(z_2)}+\lb{\a(z_1),\b(y_2)}+\lb{\b(z_1),\b(z_2)}\;,\\
    \jmath\(\lb{x_1^+,x_2^+}^+\)
    &=& \a(\lb{y_1,z_2}+\lb{z_1,y_2})+\b(\lb{z_1,z_2})\\
    &=& \a(\pb{y_1,z_2}-\s(y_1,z_2)+\pb{z_1,y_2}-\s(z_1,y_2))+\b(\pb{z_1,z_2})\\
    &=& \lb{\a(y_1),\b(z_2)}+\lb{\a(z_1),\b(y_2)}+\lb{\b(z_1),\b(z_2)}\\
    && -\a(\s(y_1,z_2)+s(z_1,y_2))\;,\\
  \s^+(x_1^+,x_2^+)\#1&=&\s(y_1,z_2)\#1+\s(z_1,y_2)\#1=\a(\s(y_1,z_2)+\a(\s(z_1,y_2))\;.
\end{eqnarray*}
This shows that
\begin{equation*}
  \lb{\jmath(x_1^+),\jmath(x_2^+)}=\gamma(\lb{x_1^+,x_2^+}^+)+\s^+(x_1^+,x_2^+)(1\#1)\;,
\end{equation*}%
for all $x_1^+,x_2^+\in\fg^+$. By the universal property of the modified Lie enveloping algebra, there
exists a (unique) algebra morphism $\gamma^{-1}$ making the following diagram commutative:
\begin{equation*}
  \begin{tikzpicture}
  \matrix (m) [matrix of math nodes,row sep=3em,column sep=4em,minimum width=2em]
  {
     \Sym(\fg)\#U_{Lie}(\fg) \\
     \fg^+ &U_{Lie,\s^+}(\fg^+)  \\};
  \path[-stealth]
    (m-2-1) edge node [left] {$\jmath$} (m-1-1)
            edge node [below] {$\iota$} (m-2-2)
    (m-2-2) edge [dashed,->] node [right] {$\gamma^{-1}$} (m-1-1);
  \end{tikzpicture}
\end{equation*}
On generators of these algebras, one checks that $\gamma$ and $\gamma^{-1}$ are inverse to each other, showing that
$\gamma$ is an algebra isomorphism (with inverse $\gamma^{-1}$).
\end{proof}


\begin{remark}\label{rem:LP_filt}
Being a (modified) Lie enveloping algebra, $U_{Lie,\s^+}(\fg^+)$ has a natural filtration, where every element of
$\fg^+$, viewed as an element of $U_{Lie,\s^+}(\fg^+)$, has filtered degree 1. As a Poisson enveloping algebra, it
also has a natural filtration; in the latter filtration, all elements of $\a'(\fg)$ have degree 0 and all elements
of $\b'(\fg)$ have degree 1; said differently, for this filtration, every element of $\fg^0\subset\fg^+$, viewed as
an element of $U_{Lie,\s^+}(\fg^+)$, has filtered degree 0, while every element of $\fg\subset\fg^+$, viewed as an
element of $U_{Lie,\s^+}(\fg^+)$, has filtered degree 1. This issue has important consequences, as we will see
when discussing the PBW theorem (see Remark \ref{rem:PBW_LP} below).
\end{remark}
  
\begin{example}
  Let $(V,\omega)$ be a symplectic vector space of dimension $2n$ over a field $\bbF$. There exists a symplectic
  basis $(X_1, Y_1,\dots, X_n, Y_n)$ of $V$ such that $\omega(X_i,Y_j)=\delta_{i,j}$ and
  $\omega(X_i,X_j)=\omega(Y_i,Y_j)=0$, for all $i,j$. Viewing $V$ as a trivial Lie algebra, $\omega$ is a 2-cocycle
  in the trivial Lie algebra cohomology of $V$. The symmetric algebra $\Sym(V)\simeq \bbF[X_1,Y_1,\dots, X_n,Y_n]$
  is a modified Lie-Poisson algebra, whose Poisson bracket is given by $\{X_i,Y_i\}_\omega=1$ ($1\leqs i\leqs n$)
  and all other brackets between basis elements are zero.  The double, $V^+=V^0\oplus V$, has as basis
  $(p_i,q_i)_{1\leqs i\leqs 2n}$ with $q_i=X_i^0$, $q_{n+i}=Y_i^0$, $p_i=-Y_i$, $p_{n+i}=X_i$ for $1\leqs i\leqs
  2n$. It follows that the Poisson enveloping algebra ${U}_{Lie,\omega}(V^+)$ of $(\Sym(V),\PB_\omega)$ is the Weyl
  algebra $A(V)\simeq A_{4n}(\bbF)$.
\end{example}

\subsubsection{The Poisson enveloping algebra of a polynomial Poisson algebra as a quotient of a smash product 
algebra}
Let $M$ be an $R$-module and suppose that its symmetric algebra $\Sym(M)$ is equipped with a Poisson bracket $\PB$,
making it into a polynomial Poisson algebra (see Example \ref{ex:poisson_poly}). We give in this paragraph a
construction of its Poisson enveloping algebra, using the smash product algebra $\Sym(M)\# T(M)$, constructed in
Example \ref{exa:smash_poly}. Let $\psi_M:\Sym(M)\to \Sym(M)\# T(M)$ denote the unique ($R$-linear) derivation of
$\Sym(M)$ with values in $\Sym(M)\# T(M)$, defined by $\psi_M(x):=1\#x$ for all $x\in M$. For a monomial
$\ux=x_1x_2\dots x_k\in \Sym(M)$,
\begin{equation}\label{eq:psi_def}
  \psi_M(\ux)=\sum_{i=1}^k\omiet xi\# x_i\;.
\end{equation}%
Notice that $\psi_M$ actually takes values in $\Sym(M)\otimes M$; it will sometimes be convenient to view $\psi_M$
as a map $\Sym(M)\to \Sym(M)\otimes M$, but we will always use the same notation $\psi_M$, because there is no risk
of confusion.

We denote by $J_M$ the two-sided ideal of $\Sym(M)\# T(M)$, generated by all elements $1\#[x,y]_\otimes-
\psi_M(\pb{x,y})$, where $x$ and $y$ both run through $M$, and where $[x,y]_\otimes:=x\otimes y-y\otimes x$. Let
$\pi_M:\Sym(M)\# T(M)\to \Sym(M)\# T(M)/J_M$ denote the canonical surjection and let $\a$ and $\b$ denote the
maps, defined by
\begin{equation*}
  \begin{array}{lcccl}
    \a&:&\Sym(M)&\to&\Sym(M)\# T(M)/J_M\\
    & & a&\mapsto&\pi_M(a\#1)\;,
  \end{array}
\end{equation*}
\begin{equation*}
  \begin{array}{lcccl}
    \b&:&\Sym(M)&\to&\Sym(M)\# T(M)/J_M\\
    & & a&\mapsto&\pi_M(\psi_M(a))\;.
  \end{array}
\end{equation*}
\begin{thm}\label{thm:envel_polynomial}
  $(\Sym(M)\# T(M)/J_M,\a,\b)$ is a Poisson enveloping algebra of the polynomial Poisson algebra~$\Sym(M)$.
\end{thm}
\begin{proof}
The verification of  items (1) -- (4) in Definition \ref{def:PEA} is very similar to the verification in the proof
of Theorem \ref{thm:envel_lie-poisson}, so we skip it here. 
We prove the universal property of $\Sym(M)\# T(M)/J_M$. Let $U'$ be any algebra and suppose that we are given
any algebra (resp.\ Lie algebra) morphisms $\a':\Sym(M)\to U'$ and $\b':\Sym(M)\to U'_L$, satisfying properties
(3') and~(4') (of Definition~\ref{def:PEA}). We show that there is a unique algebra morphism
$\gamma:\Sym(M)\#T(M)/J_M\to U'$ such that $\gamma\circ\a=\a'$ and $\gamma\circ\b=\b'$. As in the case of the
proof of Theorem \ref{thm:envel_lie-poisson}, the fact that every element of $\Sym(M)\# T(M)/J_M$ can be written as
a finite sum, where every term is the product of an element of $\Im(\a)$ with elements of $\Im(\b)$, implies
the uniqueness of the morphism~$\gamma$, if it exists; moreover, it leads to the following formula:
\begin{eqnarray*}
  \gamma(\pi_M(a\#(x_1\otimes x_2\otimes\cdots \otimes x_k)))=\a'(a).\b'(x_1).\b'(x_2)\dots \b'(x_k)\;.
\end{eqnarray*}
We need to prove that $\gamma$ is well-defined by this formula and that it is a morphism of algebras. To do this,
we first define a map $\gamma':\Sym(M)\#T(M)\to U'$ by setting
\begin{eqnarray*}
  \gamma'(a\#(x_1\otimes x_2\otimes\cdots \otimes x_k)):=\a'(a).\b'(x_1).\b'(x_2)\dots \b'(x_k)\;.
\end{eqnarray*}
The verification that $\gamma'$ is a morphism of algebras is exactly the same as the verification given in the
proof of Theorem \ref{thm:envel_lie-poisson} that $\gamma$ is a morphism of algebras.  It follows that, in order to
show that $\gamma$ is well-defined, it suffices to show that $\gamma'$ vanishes on elements of the form
$(1\#x)\odot(1\#y)-(1\#y)\odot(1\#x)- \psi_M(\pb{x,y})$, where $x,y\in M$:
\begin{eqnarray*}
\gamma'(1\#\lb{x,y}_\otimes- \psi_M(\pb{x,y}))&=&\b'(x).\b'(y)-\b'(y).\b'(x)-\gamma'(\psi_M(\pb{x,y}))\\
  &=&  \b'(\pb{x,y})-\gamma'(\psi_M(\pb{x,y}))\;.
\end{eqnarray*}
In order to show that the latter expression is zero, we show that $\b'(\ux)-\gamma'(\psi_M(\ux))=0$ for any
monomial $\ux\in M$. Since $\b'$ is an $\a'$-derivation, we have
\begin{eqnarray*}%
  \b'(\ux)&=&\sum_i\a'(\omiet x i).\b'(x_i)=\sum_i\gamma'(\omiet xi\#x_i)=\gamma'(\psi_M(\ux))\;.
\end{eqnarray*}%
Since $\gamma'$ vanishes on the ideal $J_M$ there exists a unique algebra morphism $\gamma$, such that
$\gamma\circ\pi_M=\gamma'$. In particular, $\gamma$ is defined by the above formula and satisfies
$\gamma\circ\a=\a'$ and $\gamma\circ\b=\b'$.
\end{proof}

\subsection{The Poisson enveloping algebra of a general Poisson algebra}\label{general_Uoisson}
We consider in this subsection the construction of the Poisson enveloping algebra of a general Poisson algebra and
derive from it the natural filtration of the Poisson enveloping algebra.

\subsubsection{The Poisson enveloping algebra of a general Poisson algebra as a quotient of a smash product
  algebra} The idea of the construction that we give is due to Huebschmann \cite{huebsch}, who gives an alternative
construction to Rinehart's construction of the enveloping algebra of a Lie-Rinehart algebra in terms of
Massey-Peterson algebras. We give the construction and only sketch the proof, because it is very similar to the
proof of Theorem \ref{thm:envel_polynomial}.

We start from the smash product algebra $\cA\# U_{Lie}(\cA)$, which was constructed in Example
\ref{ex:smash_general}.  Let $K$ denote the two-sided ideal of $\cA\# U_{Lie}(\cA)$ generated by all
elements\footnote{To be precise, if we denote by $\imath$ the canonical injection of $\cA$ in its Lie enveloping
  algebra and by $1_\cA$ the unit of $\cA$, then $K$ is the ideal generated by all elements of the form
  $a\#\imath(b)+b\#\imath(a)-1_\cA\#\imath(ab)$, for $a,b\in\cA$.}  of the form $a\#b+b\# a-1\#ab$ for
$a,b\in\cA$. Let $\pi_K$ denote the canonical surjection $\pi_K:\cA\# U_{Lie}(\cA) \rightarrow \cA\#
U_{Lie}(\cA)/K$ and consider the maps $\a$ and $\b$, defined by
\begin{equation*}
  \begin{array}{lcccl}
    \a&:&\cA&\to&\cA\# U_{Lie}(\cA)/K\\
    & & a&\mapsto&\pi_K(a\# 1)\;,
  \end{array}
\end{equation*}
\begin{equation*}
  \begin{array}{lcccl}
    \b&:&\cA&\to&\cA\# U_{Lie}(\cA)/K\\
    & & a&\mapsto&\pi_K(1\# a)\;.
  \end{array}
\end{equation*}
\begin{thm}\label{thm:envel_general}
  $(\cA\# U_{Lie}(\cA)/K,\a,\b)$ is a Poisson enveloping algebra of~$\cA$.
\end{thm}
\begin{proof}
The verification of properties (1) - (4) of Definition \ref{def:PEA} is not quite the same as the proof of these
properties in Theorem \ref{thm:envel_lie-poisson}, but does not pose any real difficulty. For example, (4) is now a
consequence of the definition of the ideal $K$; also, the verification of (2) is now even quicker, because it
amounts to the equality of the right hand sides of the following two formulas, valid for $a_1,a_2\in\cA$:
\begin{eqnarray*}
  \b(\pb{a_1,a_2})&=&\pi_K(1\#\pb{a_1,a_2})=\pi_K(1\#(a_1.a_2-a_2.a_1))\;,\\  
  \lb{\b(a_1),\b(a_2)}&=&\pi_K((1\#a_1)\odot(1\#a_2)-(1\#a_2)\odot(1\#a_1))\;.
\end{eqnarray*}%
The proof that $\cA\# U_{Lie}(\cA)/K$ satisfies the universal property is essentially the same as the proof which
we gave of Theorem \ref{thm:envel_lie-poisson}.
%
%
\end{proof}

\subsubsection{The filtration of the Poisson enveloping algebra}\label{ssec:filtration} 
One immediate consequence of the construction in the previous paragraph is that the Poisson enveloping algebra
$U(\cA)$ of any Poisson algebra $\cA$ is generated, as an $R$-algebra, by the images of the maps $\a$ and
$\b$. For $k\in\bbN$, we denote by $U_k(\cA)$ the $\cA$-submodule of $U(\cA)$, generated by all products of
at most $k$ elements of $\b(\cA)$.

\begin{prop}\label{lma:UP_filtered}
  Let $\cA$ be any Poisson algebra. Its Poisson enveloping algebra is a filtered $R$-algebra,
  $U(\cA)=\bigcup\limits_{i\in\bbN} U_i(\cA)$, where the filtration is given by $\cA$-submodules. Moreover, this
  filtration coincides with the filtration which is induced by the canonical filtration of $U_{Lie}(\cA)$ (taking
  on the first component $\cA$ of $\cA\# U_{Lie}(\cA)$ the trivial filtration).
\end{prop}
\begin{proof}
It follows from Theorem \ref{thm:envel_general} that $U(\cA)=\bigcup\limits_{i\in\bbN} U_i(\cA)$.  A key property
is that the elements in the images of $\a$ and $\b$ commute, modulo elements in the image of $\a$. Indeed,
according to item (3) in Definition~\ref{def:PEA} we have that $[\a(a_1),\b(a_2)]=\a(\pb{a_1,a_2})$ for any
$a_1,a_2\in\cA$. The property implies on the one hand that $U_kU_\ell\subset U_{k+\ell}$ for all
$k,\ell\in\bbN$. On the other hand, it implies that $U_k(\cA)=\pi_K(\cA\otimes U_{Lie,k}(\cA))$, since
$U_{Lie,k}(\cA)$ is by definition the $R$-module generated by products of at most $k$ elements of~$\cA$.
\end{proof}

\begin{prop}\label{prp:U_functor}
  Let $\cA$ and $\cB$ be Poisson algebras with Poisson enveloping algebras $(U(\cA),\a_\cA,\b_\cA)$ and
  $(U(\cB),\a_\cB,\b_\cB)$. For every morphism of Poisson algebras $f:\cA\to\cB$, there exists a unique morphism
  $U(f):U(\cA)\to U(\cB)$ of filtered algebras, making the following diagram commutative:
  \begin{center}
    \begin{tikzcd}
      \cA\arrow{d}[swap]{f}\arrow{r}{\a_\cA,\b_\cA}&U(\cA)\arrow{d}{U(f)}\\
      \cB\arrow{r}{\a_\cB,\b_\cB}&U(\cB)
    \end{tikzcd}
  \end{center}
\end{prop}
\begin{proof}
Uniqueness of the map $U(f)$ is clear, because $U(\cA)$ is generated by the images of $\a_\cA$ and $\b_\cA$. For
its construction, consider the maps $\a_\cB\circ f,\b_\cB\circ f:\cA\to U(\cB)$. They are algebra, resp.\ Lie
algebra morphisms satisfying the conditions (3') and (4') in Definition \ref{def:PEA} because $\a_\cB$ and $\b_\cB$
have these properties. Thus, the universal property of $U(\cA)$ yields the unique algebra morphism $U(f)$ which
completes the diagram into a commutative one. It is explicitly given by
$U(f)(a\b_\cA(a_1).\b_\cA(a_2)\dots\b_\cA(a_k))=f(a)\b_\cB(f(a_1)).\b_\cB(f(a_2))\dots \b_\cB(f(a_k))$, where
$a,a_1,a_2,\dots,a_k\in\cA$. From this formula it is clear that $U(f)$ is filtered.
\end{proof}
The proposition implies that there is a covariant functor $U$ between the category of Poisson algebras over $R$ and
the category of filtered algebras over~$R$, which assigns to each Poisson algebra $\cA$ the Poisson enveloping
algebra~$U(\cA)$, given by Theorem \ref{thm:envel_general} and to each morphism of Poisson algebras $f:\cA\to\cB$
the induced morphism $U(f):U(\cA)\to U(\cB)$, given by Proposition~\ref{prp:U_functor}.


\subsection{The Poisson enveloping algebra of a quotient of a general Poisson algebra}\label{par:UP_quotient}
Suppose that $\cA$ is any Poisson algebra and that $I$ is a Poisson ideal of $\cA$.  We give in this subsection a
description of the Poisson enveloping algebra of the Poisson algebra $\cB:=\cA/I$ in terms of the Poisson
enveloping algebra $(U(\cA),\a_\cA,\b_\cA)$ of~$\cA$. To do this, we first construct a new algebra out of $\cB$ and
$U(\cA)$. Using the canonical surjection $\pi:\cA\to\cB$ and the commutativity of $\cA$ and $\cB$, we make $\cB$
into a symmetric $\cA$-module by setting $b\cdot a=a\cdot b:=\pi(a)b$, for $a\in \cA$ and $b\in \cB$. Consider the
$\cB$-module $\cB\otimes_\cA U(\cA)$. If we denote the unit of $\cB$ by~$1_\cB$, then in $\cB\otimes_\cA U(\cA)$ we
have the equality $\pi(a)\otimes u=1_\cB\otimes au=1_\cB\otimes\a_\cA(a).u$, valid for all $a\in\cA$ and $u\in
U(\cA)$.
\begin{prop}\label{prop:nul_on_ideal}
  The $\cB$-module $\cB\otimes_\cA U(\cA)$ is a unitary algebra over $R$ with the product defined for
  $\pi(a_i)\otimes u_i\in\cB\otimes_\cA U(\cA)$, $i=1,2$, by
  \begin{equation}\label{for:new_product_def}
    (\pi(a_1)\otimes u_1)\cdot( \pi(a_2)\otimes u_2):= 1_\cB\otimes (a_1u_1).(a_2u_2)=1_\cB\otimes 
    \alpha_\cA(a_1).u_1.\alpha_\cA(a_2).u_2\;.
  \end{equation}%
\end{prop}
\begin{proof}
If the product is well-defined then it is clear that it is associative and has $1_\cB\otimes 1_{U(\cA)}$ as unit.
To prove that the product is well-defined it is sufficient to show that if $\jmath\in I$ and $u,v\in U(\cA)$ then,
in the $\cB$-module $\cB\otimes_\cA U(\cA)$, one has $1_\cB\otimes u.\a_\cA(\jmath).v=0$.  In view of Proposition
\ref{lma:UP_filtered}, it suffices to show that $1_\cB\otimes u.\a_{\cA}(\jmath).v=0$ for all $u$ of the form
$\b_\cA(a_1).\b_\cA(a_2)\dots\b_\cA(a_k)$, where $k\in\bbN$ and all $a_i$ belong to $\cA$. We do this by recursion
on $k$. Since $1_\cB\otimes \a_\cA(\jmath).v=1_\cB\otimes \jmath v=\pi(\jmath)\otimes v=0$ for any $\jmath\in I$,
the result is clair for $k=0$. Let us assume it to be true up to order $k-1$. Using successively the relation
$\b_\cA(a_i).\a_\cA(\jmath)=\a_\cA(\jmath).\b_\cA(a_i)+\a_\cA(\pb{a_i,\jmath})$ for $i=k,k-1,\dots,1$ allows us to
permute $\a_\cA(\jmath)$ with all $\b_\cA(a_i)$, because $\pb{a_i,\jmath}\in I$ (recall that $I$ is a Poisson
ideal), and so each one of the correction terms $1_\cB\otimes \b_\cA(a_1).\b_\cA(a_2)\dots
.\b_\cA(a_{i-1}).\a_\cA(\pb{a_i,\jmath}).u'.v$ is zero, in view of the recursion hypothesis. It follows that
\begin{equation*}
  1_\cB\otimes u.\a_\cA(\jmath).v=1_\cB\otimes \jmath u.v=\pi(\jmath)\otimes u.v=0\;.
\end{equation*}%
\end{proof}
Let
$I_\cB$ denote the (two-sided) ideal of $\cB\otimes_\cA U(\cA)$, generated by $1_\cB\otimes\b_\cA(I)$ and let
$\pi_\cB: \cB\otimes_\cA U(\cA)\to \cB\otimes_\cA U(\cA)/I_\cB$ denote the canonical surjection. Consider the maps
$\a$ and $\b$, defined by
\begin{equation}\label{eq:alpha_def}
  \begin{array}{lcccl}
    \a&:&\cB&\to&\cB\otimes_\cA U(\cA)/I_\cB\\
    & & \pi(a)&\mapsto&\pi_\cB(1_\cB\otimes \a_\cA(a))\;,
  \end{array}
\end{equation}
\begin{equation}\label{eq:beta_def}
  \begin{array}{lcccl}
    \b&:&\cB&\to&\cB\otimes_\cA U(\cA)/I_\cB\\
    & & \pi(a)&\mapsto&\pi_\cB(1_\cB\otimes\b_\cA(a))\;.
  \end{array}
\end{equation}
Notice that both maps all well-defined: the first one in view of the proof of proposition \ref{prop:nul_on_ideal}
and the second one in view of the definition of the ideal $I_\cB$.
\begin{thm}\label{thm:envel_general_quotient}
	Let $\cA$ be any Poisson algebra, $I$ a Poisson ideal of $\cA$, and $(U(\cA),\a_\cA,\b_\cA)$ a Poisson enveloping algebra of $\cA$. A Poisson enveloping algebra of the quotient Poisson algebra
	$\cB=\cA/I$ is $(\cB\otimes_\cA U(\cA)/I_\cB,\a,\b)$.
\end{thm}
\begin{proof}
The proof that $\a$ and $\b$ satisfy the properties (1) -- (4) in Definition~\ref{def:PEA} is an
immediate consequence of the fact that $\a_\cA$ and $\b_\cA$ satisfy these properties, in combination with the
following three formulas, which are a direct consequence of definition (\ref{for:new_product_def}): for any
$a_1,a_2\in\cA$ and $u_1,u_2\in U(\cA)$, 
\begin{eqnarray*}
  (\pi(a_1)\otimes 1_{U(\cA)})\cdot( \pi(a_2)\otimes 1_{U(\cA)})&=&\pi(a_1a_2)\otimes 1_{U(\cA)}\;,\\
  (1_\cB\otimes u_1)\cdot( 1_\cB\otimes u_2)&=&1_\cB\otimes u_1.u_2\;,\\
  (\pi(a_1)\otimes 1_{U(\cA)})\cdot( 1_\cB\otimes u_2)&=&\pi(a_1)\otimes u_2\;.
\end{eqnarray*}
In order to show that $U(\cB):=\cB\otimes_\cA U(\cA)/I_\cB$ satisfies the universal property, suppose that $U'$ is
any algebra and that $\a':\cB\rightarrow U',\b':\cB\rightarrow U'_L$ are algebra (resp. Lie algebra) morphisms,
satisfying (3') and (4') in Definition \ref{def:PEA}. We will prove that there exists a unique algebra morphism
$\gamma:U(\cB)\rightarrow U'$ such that $\gamma\circ\a=\a'$ and $\gamma\circ\b=\b'$. We do this by showing that the
following diagram is a commutative diagram: the above relations which $\gamma$ is ought to satisfy are equivalent
to the commutativity of the triangle $(5)$.
\begin{center}
  \begin{tikzcd}
    U(\cA)
    \arrow{ddr}{\iota}
    \arrow[bend left=20, dotted]{ddddrrrrr}{\gamma'}[swap,xshift=0ex,yshift=-2ex]{(3)}&&&&&\\
    \\
    &\cB\otimes_\cA U(\cA)
    \arrow{d}[swap]{\pi_\cB}
    \arrow[dotted]{ddrrrr}{\gamma''}[swap,xshift=-1.5ex,yshift=0.5ex]{(4)}&&&&\\
    &U(\cB)
    \arrow[dotted]{drrrr}[swap]{\gamma}&&&&\\
    \cA
    \arrow{ur}[rotate=23,anchor=center,yshift=1.3ex]{\overline{\a},\overline{\b}}[swap]{(1)}
    \arrow{uuuu}{\a_\cA,\b_\cA}[swap,xshift=4ex,yshift=-1ex]{(2)}
    \arrow{r}[swap]{\pi}&
    \cB
    \arrow{u}{\a,\b}
    \arrow{rrrr}[swap]{\a',\b'}[xshift=-4ex,yshift=1ex]{(5)}&&&&U'
  \end{tikzcd}
\end{center}
The morphisms $\a'\circ\pi$ and $\b'\circ\pi$ are algebra (resp. Lie algebra) morphisms and satisfy the same
properties (3') and (4') as $\a'$ and $\b'$, so by the universal property of $U(\cA)$ there exists a (unique)
algebra morphism $\gamma':U(\cA)\to U'$ which makes the outer diagram commute. Consider the linear map
$\iota:U(A)\rightarrow \cB\otimes_\cA U(A)$ defined for all $u\in U(A)$ by $\iota(u)=1_\cB\otimes
u$. Proposition~\ref{prop:nul_on_ideal} shows that $\iota$ is an algebra morphism. Consequently,
$\overline{\a}=\pi_\cB\circ\iota\circ \a_\cA$ and $\overline{\b}=\pi_\cB\circ\iota\circ \b_\cA$ are algebra
(resp. Lie algebra) morphisms. The definition of $\a$ and $\b$ implies that the diagrams $(1)$ and $(2)$ commute.

\noindent Consider the linear map
\begin{equation*}
  \begin{array}{lcccl}
    \gamma''&:&\cB\otimes_\cA U(\cA)&\to& U'\\
    & & b\otimes u&\mapsto&\a'(b).\gamma'(u)\;.
  \end{array}
\end{equation*}
For $a\in\cA$ we have that 
\begin{eqnarray*}
	\gamma''((a\cdot b)\otimes u)&=&\gamma''((b\cdot a)\otimes u)=\a'(b\pi(a)).\gamma'(u)\\
	&=&\a'(b).\a'(\pi(a)).\gamma'(u)=\a'(b).\gamma'(\a_\cA(a)).\gamma'(u)\\
	& =&\gamma''(b\otimes (au))\;,
\end{eqnarray*}
so that $\gamma''$ is well-defined. If we make $U'$ into a $\cB$-module upon using $\a'$, then $\gamma''$ can be
described as the unique morphism of $\cB$-modules, which sends $\iota(u)=1_\cB\otimes u$ to $\gamma'(u)$. Thus the
diagram $(3)$ commutes. Since $\gamma'$ is an algebra morphism, it follows from this description that $\gamma''$ is
also an algebra morphism. For $\jmath\in I$ we have
\begin{equation*}
  \gamma''(1_\cB\otimes \b_\cA(\jmath))=\gamma'(\b_\cA(\jmath))=\b'(\pi(\jmath))=0\;,
\end{equation*}%
so that $\gamma''$ induces an algebra morphism $\gamma:\cB\otimes_\cA U(\cA)/I_\cB\to U'$, such that the diagram
$(4)$ commutes. The commutativity of the diagrams (1) -- (4) and of the outer diagram shows that $\gamma\circ
\a\circ\pi=\a'\circ \pi$ and $\gamma\circ \b\circ\pi=\b'\circ \pi$. By surjectivity of $\pi$ we conclude that the
diagram $(5)$ commutes. It remains to be shown that the morphism $\gamma$ is unique. This follows from the fact
that $\cB\otimes_\cA U(\cA)/I_\cB$ is generated by the images of $\a$ and $\b$, which is in turn a consequence of
the fact that $U(\cA)$ is generated by the images of $\a_\cA$ and~$\b_\cA$.
\end{proof}
\begin{remark}
Let $I_P$ denote the two-sided ideal of $U(\cA)$, generated by
$\a_\cA(I)$ and $\b_\cA(I)$. It can be shown as above that $(U(\cA)/I_P,\a,\b)$ is a Poisson enveloping
algebra for $\cB=\cA/I$, where $\a$ and $\b$ are defined as the unique morphisms which make the following
diagram commutative:
\begin{center}
  \begin{tikzcd}
    \cA\arrow{r}{\a_\cA,\b_\cA}\arrow{d}[swap]{\pi}&U(\cA)\arrow{d}{}\\
    \cB\arrow{r}{\a,\b}&\frac{U(\cA)}{I_P}
  \end{tikzcd}
\end{center}
In this diagram, the vertical arrows are the canonical surjections.
\end{remark}

\section{The Poincar\'e-Birkhoff-Witt theorem}\label{sec:PBW}

\subsection{The graded algebra associated with $U(\cA)$}
Let $\cA$ be a Poisson algebra (over $R$) and let $(U(\cA),\a,\b)$ be its Poisson enveloping algebra. We recall
from Section \ref{ssec:filtration} that $U(\cA)$ has a canonical filtration,
$$
  U(\cA)=\bigcup\limits_{i\in\bbN} U_i(\cA)\;, 
$$
where $U_k(\cA)$ stands for the $\cA$-submodule of $U(\cA)$, generated by all products of at most $k$ elements of
$\b(\cA)$, where $k\in\bbN$. The graded algebra (over~$R$) associated with the filtered algebra $U(\cA)$ is given
by
\begin{equation*}
  \Gr(U(\cA))=\bigoplus_{i\in\bbN}\Gr^i(U(\cA)),\quad\hbox{where}\quad 
  \Gr^k(U(\cA)):=\frac{U_k(\cA)}{U_{k-1}(\cA)}\;.
\end{equation*}%
The homogeneous components $\Gr^k(U(\cA))$ are $\cA$-modules, just like the $\cA$-submodules $U_k(\cA)$ of $U(\cA)$
from which they are constructed. As in the case of Lie algebras, we have the following result:

\begin{prop}
  $\Gr(U(\cA))$ is a commutative $\cA$-algebra.
\end{prop}
\begin{proof}
In terms of the canonical surjections 
\begin{equation}\label{eq:rho}
  \Gr_k:U_k(\cA)\to\frac{U_k(\cA)}{U_{k-1}(\cA)}\;,
\end{equation}
the product on $\Gr(U(\cA))$ is given, for $\xi_k\in U_k(\cA)$ and $\xi_\ell\in U_\ell(\cA)$ by
$\Gr_k(\xi_k)\Gr_\ell(\xi_\ell):=\Gr_{k+\ell}(\xi_k.\xi_\ell)$. The fact that $\b$ is a Lie algebra morphism, item
(3) in Definition \ref{def:PEA} and the commutativity of $\cA$ imply respectively that
\begin{eqnarray*}
  \lb{\Gr_1(\b(a_1)),\Gr_1(\b(a_2))}
    &=&\Gr_2(\lb{\b(a_1),\b(a_2)})=\Gr_2(\b(\pb{a_1,a_2}))=0\;,\\ 
  \lb{\Gr_1(\a(a_1)),\Gr_1(\b(a_2))}
    &=&\Gr_2(\lb{\a(a_1),\b(a_2)})=\Gr_2(\a(\pb{a_1,a_2}))=0\;,\\ 
  \lb{\Gr_1(\a(a_1)),\Gr_1(\a(a_2))}&=&0\;,
\end{eqnarray*}
for all $a_1,a_2\in\cA$. It follows that the product on $\Gr(U(\cA))$ is $\cA$-bilinear and commutative.
\end{proof}
\subsection{The PBW map}\label{par:PBW-map}
As in the previous section, let $\cA$ be a Poisson algebra and $(U(\cA),\a,\b)$ its Poisson enveloping
algebra. Recall that we denote by $\Omega(\cA)$ the $\cA$-module of K\"ahler differentials of $\cA$. Recall also
that (4) in Definition \ref{def:PEA} says that $\b$ is an $\a$-derivation of $\cA$ with values in $U(\cA)$. The
universal property of $\Omega(\cA)$ leads to an $\cA$-linear map, defined by
\begin{equation}\label{eq:psi}
  \begin{array}{lcccl}
    \psi&:&\Omega(\cA)&\to&U(\cA)\\
        & &\diff a &\mapsto&\b(a)\;.
  \end{array}
\end{equation}
Notice that $\psi$ actually takes values in $U_1(\cA)$. Let $\Psi:\Omega(\cA)\to\Gr(U(\cA))$ be the induced map,
which is a morphism of graded $\cA$-modules with values in a commutative $\cA$-algebra. By the universal property
of the symmetric algebra $\Sym_\cA(\Omega(\cA))$ of $\Omega(\cA)$ we get a morphism of graded $\cA$-algebras
$\Sym_\cA(\Omega(\cA))\to \Gr(U(\cA))$. It is called the \emph{Poincar\'e-Birkhoff-Witt map}, or \emph{PBW map} for
short, and is explicitly given by
\begin{equation}\label{eq:PBW_map}
  \begin{array}{lcccl}
    \PBW_\cA&:&\Sym_\cA(\Omega(\cA))&\to&\Gr(U(\cA))\\
        & &a\diff a_1 \diff a_2\dots \diff a_k &\mapsto&\Gr_k(a\b(a_1).\b(a_2)\dots\b(a_k))\;,
  \end{array}
\end{equation}
where $a,a_1,\dots,a_k\in\cA$.  The latter image can also be written as the product
$a\Gr_1(\b(a_1)).\Gr_1(\b(a_2))\dots\Gr_1(\b(a_k))$.  It is clear from (\ref{eq:PBW_map}) and
Proposition~\ref{lma:UP_filtered} that the PBW map is surjective. Also, $\PBW_\cA$ is a map of \emph{graded}
$\cA$-algebras, because $\Gr_1({\b(a)})$ is homogeneous of degree 1 in $\Gr(U(\cA))$ for any $a\in\cA$.

\begin{defn}%
  A Poisson algebra $\cA$ \emph{satisfies the PBW theorem} if the graded map $\PBW_\cA:\Sym_\cA(\Omega(\cA))
  \to\Gr(U(\cA))$ is injective, hence is an isomorphism of graded $\cA$-algebras.
\end{defn}
At this moment we do not know of any Poisson algebra which does \emph{not} satisfy the PBW theorem. We give a few
examples here and elaborate on some other examples in the subsections that follow.
\begin{example}\label{ex:smooth}
Any smooth Poisson algebra $\cA$ over a field $\bbF$ satisfies the PBW theorem: the pair $(\cA,\Omega(\cA))$ is a
Lie-Rinehart algebra (see \cite{huebsch}) and Rinehart shows in \cite{rinehart} that the PBW theorem holds for
Lie-Rinehart algebras $(A,L)$ under the condition that $L$ is a projective $A$-module; in our case, $\Omega(\cA)$
is a projective $\cA$-module because $\cA$ is assumed to be a smooth algebra over a field.
\end{example}
\begin{example}
Let $\cA$ be any algebra which we make into a Poisson algebra by adding the zero Poisson bracket. We have shown in
Proposition \ref{prop:null} that $\Sym_\cA(\Omega(\cA))$ is a Poisson enveloping algebra of $\cA$. In this case,
the enveloping algebra is already graded (with grading coming indeed from the canonical filtration of the Poisson
enveloping algebra) and the map $\PBW_\cA$ is just the identity map. In particular, $\cA$ satisfies the PBW
theorem. Notice that this example is not a particular case of the previous one: here $\cA$ can be any algebra,
smooth or singular.
\end{example}

\subsection{The PBW theorem for modified Lie-Poisson algebras}
We show in this subsection that the PBW theorem holds for any (modified) Lie-Poisson algebra when the base ring $R$
is a field $\bbF$ (see Example \ref{ex:affine}). The result is not new, as such an algebra is obviously smooth, so
it is covered by Example \ref{ex:smooth}, but our proof which is specific to the Lie-Poisson case, has some extra
flavors, such as being more direct, more explicit, and it prepares for the singular case, which we will study in
the next subsection.

%
%
%
\begin{thm}\label{thm:PBW_LP}
  Let $\fg$ be a Lie algebra over $\bbF$ and let $\sigma$ be a 2-cocycle in the trivial Lie algebra cohomology of
  $\fg$. The modified Lie-Poisson algebra $(\Sym(\fg),\PB_\sigma)$ satisfies the PBW theorem.
\end{thm}
\begin{proof}
In order to simplify the notation, we denote throughout this proof $\Sym(\fg)$ by $\cA$, adding a subscript
$\sigma$ to $\cA$ when the Poisson structure $\PB_\s$ is relevant.  We need to show that the PBW map
\begin{equation*}
  PBW_{\cA_\sigma}:\Sym_{\cA}(\Omega(\cA))\to\Gr(U(\cA_\sigma))
\end{equation*}%
is an isomorphism of graded $\cA$-algebras. It follows from the following chain of isomorphisms of graded
$\cA$-algebras, each of which will be detailled below:
\begin{eqnarray*}
 \Gr(U(\cA_\s))&\stackrel{(1)}=&\Gr(\cA\#U_{Lie}(\fg))\stackrel{(2)}=\Gr(\cA\otimes U_{Lie}(\fg)) \\ 
   &\stackrel{(3)}\simeq&\Gr(U_{Lie}(\cA\otimes \fg))\stackrel{(4)}\simeq\Sym_{\cA}(\cA\otimes\fg)\\ 
   &\stackrel{(5)}\simeq&\Sym_{\cA}(\Omega(\cA))\;.
\end{eqnarray*}
We have shown in Section \ref{par:envel_lie-poisson} that the Poisson enveloping algebra of $\cA_\sigma$ is given
by
\begin{equation}\label{for:U_Lie_rappel}
  U(\cA_\sigma)=\Sym(\fg)\#U_{Lie}(\fg)=\cA\#U_{Lie}(\fg)\;.
\end{equation}%
This leads to the proof of (1). 

Recall that the smash product algebra $\cA\#U_{Lie}(\fg)$ is the tensor product $\cA\otimes U_{Lie}(\fg)$, with a
special product (dictated by $\sigma$ and the bracket of $\fg$). Also, the filtration of $\cA\#U_{Lie}(\fg)$ is
induced by the filtration of $U_{Lie}(\fg)$, just like the filtration on $\cA\otimes U_{Lie}(\fg)$. Therefore, the
proof of (2) amounts to showing that the product of two elements in $\cA\#U_{Lie}(\fg)$ is their tensor product,
modulo terms of lower degree. To do this, let $a_1,a_2\in \cA$ and let $x,x_1,\dots,x_k\in\fg$. According to
(\ref{for:smash_lie-poisson}),
\begin{equation}\label{for:smash_lie-poisson_2}
  (a_1\# x)\odot(a_2\# x_1.x_2\dots x_k) =a_1a_2\#x.x_1\dots x_k+a_1\pb{x,a_2}_\sigma\#x_1.x_2\dots x_k\;,  
\end{equation}%
where the first term is just the tensor product of the two arguments and the second term belongs to $(\cA\otimes
U_{Lie}(\fg))_k=\cA\otimes U_{Lie,k}(\fg)$. This shows the claim for $(a_1\#u)\odot(a_2\#v)$ with $u\in\fg\subset
U_{Lie}(\fg)$; (2) follows from it by writing a general homogeneous element $a\#u$ as the product of elements of
the form $a_1\#x$ with $a_1\in\cA$ and $x\in\fg$ and repeatedly using (\ref{for:smash_lie-poisson_2}).  

In order to prove~(3), first notice that $\cA\otimes \fg$ is a Lie algebra over $\cA$ by \emph{extension of
  scalars} (see \cite{bourbakiLA1}), namely a Lie bracket on $\cA\otimes \fg$ is given for $a_1,a_2\in\cA$ and
$x_1,x_2\in\fg$ by $\lb{a_1\otimes x_1,a_2\otimes x_2}:=a_1a_2\otimes\lb{x_1,x_2}$. The natural inclusion $\fg\to
U_{Lie}(\fg)$ leads to an inclusion $\cA\otimes\fg\to \cA\otimes U_{Lie}(\fg)$, hence (by the universal property of
the enveloping algebra) to an algebra morphism $U_{Lie}(\cA\otimes\fg)\to \cA\otimes U_{Lie}(\fg)$. Its inverse is
the unique morphism of algebras which sends $1_\cA\otimes x\in \cA\otimes U_{Lie}(\fg)$ to $1_\cA\otimes x\in
U_{Lie}(\cA\otimes \fg)$. Since the isomorphism respects the (natural) filtrations on $U_{Lie}(\cA\otimes\fg)$ and
$\cA\otimes U_{Lie}(\fg)$, this shows (3).

The isomorphism~(4) is just the classical PBW theorem for the $\cA$-Lie algebra $\cA\otimes \fg$: the PBW theorem
holds for all Lie algebras over $\cA$ because $\cA$ contains a field (the field $\bbF$) (see \cite{cohn}). Finally,
for any $\bbF$-vector space~$V$, the module $\Omega(\Sym(V))$ is a free $\Sym(V)$-module, to wit,
$\Omega(\Sym(V))\simeq \Sym(V)\otimes V$ (see \cite[Ch.\ 16]{eisenbud}). Applied to $V=\fg$ we get (5).
\end{proof}
For simplicity, we have assumed in Theorem \ref{thm:PBW_LP} that $\fg$ is a Lie algebra over a field. The
assumption was used to assert that the $\cA$-Lie algebra $\cA\otimes \fg$ satisfies the PBW theorem (for Lie
algebras!): the above proof works under the latter, more general, assumption.

\begin{remark}\label{rem:PBW_LP}
  As we have seen in Section \ref{par:envel_lie-poisson_2}, the Poisson enveloping algebra of
  $(\Sym(\fg),\PB_\sigma)$  can also be described as a (modified) Lie enveloping algebra and one might be tempted
  to use the PBW theorem for (modified) Lie algebras to show that $(\Sym(\fg),\PB_\sigma)$ satisfies the PBW
  theorem. However, as we pointed out in Remark \ref{rem:LP_filt}, the filtration of $U_{Lie,\s^+}(\fg^+)$ is
  different whether we consider it as a Lie enveloping algebra or as a Poisson enveloping algebra. This means that
  the associated graded algebras and the associated PBW maps are different, and so the classical PBW theorem cannot
  be applied directly to give a quick proof of Theorem \ref{thm:PBW_LP}.
\end{remark}  

\subsection{The PBW theorem for some singular Poisson algebras}
The purpose of this subsection is to show that if $I$ is a Poisson ideal of a smooth Poisson algebra $\cA$, which
is generated (as an ideal) by a single element and such that $\cA/I$ is an integral domain, then the Poisson
algebra $\cB:=\cA/I$ satisfies the PBW theorem. We denote as before by $\pi:\cA\to\cB=\cA/I$ the canonical
surjection and we write $U(\cA)$ for the Poisson enveloping algebra of~$\cA$, with accompanying maps denoted by
$\a_\cA$ and $\b_\cA$. We recall from Section~\ref{par:UP_quotient} that $(\cB\otimes_\cA
U(\cA)/I_\cB,\a_\cB,\b_\cB)$ is a Poisson enveloping algebra of $\cB$, where the product on $\cB\otimes_\cA U(\cA)$
is given by (\ref{for:new_product_def}), the ideal $I_\cB$ is generated by $1_\cB\otimes\b_\cA(I)$, and the
morphisms $\a_\cB$ and $\b_\cB$ are given by (\ref{eq:alpha_def}) and (\ref{eq:beta_def}).

\begin{thm}\label{thm:PBW_singular}
  The Poisson algebra $\cB=\cA/I$ satisfies the PBW theorem.
\end{thm}
\begin{proof}
We first outline the proof. Consider the following diagram of graded $\cB$-algebras:
\begin{center}
  \begin{tikzcd}
    0\arrow{r}&\Gr(I_\cB)\arrow{r}{\imath_U}&\cB\otimes_\cA\Gr(U(\cA))\arrow{r}{\pi_U}&\Gr(U(\cB))\arrow{r}{}&0\\
    0\arrow{r}&(1_\cB\otimes\diff I)\arrow{r}[swap]{\imath_S}\arrow{u}{\theta}
        &\cB\otimes_\cA\Sym_\cA(\Omega(\cA))\arrow{u}[swap]{\Id_\cB\otimes\PBW_\cA}{\simeq}\arrow{r}[swap]{\pi_S}
        &\Sym_\cB(\Omega(\cB))\arrow[twoheadrightarrow]{u}[swap]{\PBW_\cB}\arrow{r}{}&0
  \end{tikzcd}
\end{center}
The construction of the different arrows will be discussed below and we will show that the diagram has exact rows
and is commutative. We need to show that the rightmost arrow $\PBW_\cB$ is injective. The middle arrow is an
isomorphism (because $\cA$ is smooth, so the PBW theorem holds for~$\cA$). We will show that the leftmost arrow
$\theta$ is surjective. By a simple diagram chase, this implies that $\PBW_\cB$ is injective: if
$Z\in\Sym_\cB(\Omega(\cB))$ is in the kernel of $\PBW_\cB$ then there exist elements
$Y\in\cB\otimes_\cA\Sym_\cA(\Omega(\cA))$ and $X'\in\Gr(I_\cB)$, such that $\pi_S(Y)=Z$ and
$\Id\otimes\PBW_\cA(Y)=\imath_U(X')$. By surjectivity of $\theta$, there exists $X\in (1_\cB\otimes \diff I)$ such
that $\theta(X)=X'$. By the commutativity and exactness properties of the diagram, we can conclude that
$Z=\pi_S(Y)=\pi_S(\imath_S(X))=0$, as was to be shown.

We now get to the details of the proof. 

\underline{Step 1:} \emph{Exactness of the bottom line}. The conormal sequence for Kähler differentials (see
\cite[Proposition 16.3]{eisenbud}), applied to the canonical surjection $\pi:\cA\to\cB$, is the exact sequence of
$\cB$-modules, given by 
\begin{equation*}
  I/I^2\to \cB\otimes_\cA\Omega(\cA)\to\Omega(\cB)\to0\;,
\end{equation*}%
where the first map sends $j\mod I^2\in I/I^2$ to $\diff j$ and the second map sends $\pi(a_1)\otimes a_2\diff a_3$
to $\pi(a_1a_2)\diff \pi(a_3)$. Since the image of the first map is the $\cB$-submodule $\langle 1_\cB\otimes\diff
I\rangle$ of $\cB\otimes_\cA\Omega(\cA)$, generated by $1_\cB\otimes\diff I$, we have the following short exact
sequence of $\cB$-modules:
\begin{equation*}
  0\to \langle 1_\cB\otimes\diff I\rangle\to \cB\otimes_\cA\Omega(\cA)\to\Omega(\cB)\to0\;.
\end{equation*}%
Applying the $\Sym$ functor, we get according to \cite[\P 6.2, Proposition 4]{bourbaki_algebra} the following
short exact sequence
\begin{equation}\label{eq:step1_ses}
  0\to (1_\cB\otimes\diff I)\to \Sym_\cB(\cB\otimes_\cA\Omega(\cA))\to\Sym_\cB\Omega(\cB)\to0\;,
\end{equation}%
where we recall that $(1_\cB\otimes\diff I)$ stands for the two-sided ideal (in this case of the $\cB$-algebra
$\Sym_\cB(\cB\otimes_\cA\Omega(\cA))$), generated by $1_\cB\otimes\diff I$, so it is a homogeneous ideal (generated
by elements of degree 1). By extension of the ring of scalars (see \cite[\P 6.4, Proposition 7]{bourbaki_algebra})
we have the following isomorphism of $\cB-$modules:
\begin{equation*}
  \Sym_\cB(\cB\otimes_\cA\Omega(\cA))\simeq \cB\otimes_\cA\Sym_\cA(\Omega(\cA))\;.
\end{equation*}%
Substitued in (\ref{eq:step1_ses}) we get the desired exactness of the bottom line. For future reference, note that
the surjection
\begin{equation*}
  \pi_S:\cB\otimes_\cA\Sym_\cA(\Omega(\cA))\to\Sym_\cB(\Omega(\cB))
\end{equation*}%
is explicitly given by 
\begin{equation*}
  \pi_S(b\otimes a\;\diff a_1\diff a_2\dots\diff a_k)=\pi(a)b\;\diff\pi(a_1)\diff\pi(a_2)\dots\diff\pi(a_k)\;,  
\end{equation*}%
where $a,a_1,\dots,a_k\in\cA$ and $b\in\cB$.

\underline{Step 2:} \emph{Exactness of the top line.} Theorem \ref{thm:envel_general_quotient} shows exactness of
the following sequence of $R$-algebras and $\cB$-modules:
\begin{center}
	\begin{tikzcd}
		0\arrow{r}&I_\cB\arrow{r}{}&\cB\otimes_\cA U(\cA)\arrow{r}{\pi_\cB}&U(\cB)\arrow{r}{}&0\;.
	\end{tikzcd}
\end{center}
The filtration $U(\cA)=\bigcup\limits_{i\in\bbN} U_i(\cA)$ by $\cA$-modules induces a filtration $\cB\otimes_\cA
U(\cA)=\bigcup\limits_{i\in\bbN} \cB\otimes_\cA U_i(\cA)$ by $\cB$-modules.  On $I_\cB$ we take the induced
filtration, i.e., $I_{\cB,k}=I_\cB\cap(\cB\otimes_\cA U_k(\cA))$, for all $k$; also, the quotient filtration on
$U(\cB)$ is the canonical filtration of $U(\cB)$ as a Poisson enveloping algebra, $\pi_\cB(\cB\otimes_\cA
U_k(\cA))= U_k(\cB)$, for all $k$. Therefore, taking the induced morphism on the graded modules and algebras, we
get the exact sequence,
\begin{center}
\begin{tikzcd}
  0\arrow{r}&\Gr(I_\cB)\arrow{r}{}&\Gr(\cB\otimes_\cA U(\cA))\arrow{r}{\Gr(\pi_\cB)}&\Gr(U(\cB))\arrow{r}{}&0\;.
\end{tikzcd}
\end{center}
Finally the graded $\cB$-algebra $\Gr(\cB\otimes_\cA U(\cA))$ is naturally isomorphic to $\cB\otimes_\cA
\Gr(U(\cA))$. Therefore we get the exact sequence of the top line of the above diagram.

\underline{Step 3:} \emph{Commutativity of the diagram.} The commuting diagrams (1) and (2) of Theorem
\ref{thm:envel_general_quotient} show that the map $\pi_{\cB}\circ \iota:U({\cA})\longrightarrow U({\cB})$
satisfies the universal property of Proposition \ref{prp:U_functor}. Uniqueness of the morphism $U(\pi)$ leads to
the equality $U(\pi)=\pi_{\cB}\circ \iota$. According to the  definition of $U(\pi)$, 
$$
  U(\pi)(a\beta_{\cA}(a_1).\beta_{\cA}(a_2) \dots \beta_{\cA}(a_k))=
  \pi(a)\beta_{\cB}(\pi(a_1)).\beta_{\cB}(\pi(a_2))\dots \beta_{\cB}(\pi(a_k))\;,
$$
while
$$
  \begin{array}{lll}
    \pi_{\cB}\circ \iota(a\beta_{\cA}(a_1).\beta_{\cA}(a_2)\dots \beta_{\cA}(a_k)) 
      & = & \pi_{\cB}(1_{\cB}\otimes a\beta_{\cA}(a_1).\beta_{\cA}(a_2)\dots \beta_{\cA}(a_k)) \\ 
      & = & \pi_{\cB}(\pi(a)\otimes\beta_{\cA}(a_1).\beta_{\cA}(a_2)\dots \beta_{\cA}(a_k)).
  \end{array} 
$$
From the equality  $U(\pi)=\pi_{\cB}\circ\iota$, we conclude that 
$$ \pi_\cB(\pi(a)\otimes \b_{\cA}(a_1).\b_{\cA}(a_2)\dots
\b_{\cA}(a_k))=\pi(a)\beta_\cB(\pi(a_1)).\beta_\cB(\pi(a_2))\dots\beta_{\cB}(\pi(a_k))\;.
$$
Let us denote ${\rm gr}'_k:U_k(B)\longrightarrow U_k(B)/U_{k-1}(B)$. For $b\in B$ and $u\in U_k(A)$, the value in
$b\otimes{\rm gr}_k(u)$ of the map $\pi_U:={\rm gr}(\pi_{\cB})$ is
$$
  \pi_U(b\otimes{\rm gr}_k(u))={\rm gr}'_k(\pi_B(b\otimes u))\;.
$$
It follows that, for $Y:=b\,\otimes\, a\,\diff a_1\diff a_2\dots\diff a_k\in \cB\otimes_\cA
\Sym_\cA(\Omega(\cA))$,
$$
\begin{array}{lll}
  \pi_U({\Id}_{\cB}\otimes {\rm PBW}_{\cA}(Y))& = & \pi_U(b\otimes a{\rm gr}_k(\beta_{\cA}(a_1).
    \beta_{\cA}(a_2)\dots\beta_{\cA}(a_k))) \\
    & = & {\rm gr}'_k(\pi_{\cB}(b\otimes a\beta_{\cA}(a_1).\beta_{\cA}(a_2)\dots \beta_{\cA}(a_k))) \\
    & = & {\rm gr}'_k(b\pi(a)\beta_{\cB}(\pi(a_1)).\beta_{\cB}(\pi(a_2))\dots\beta_{\cB}(\pi(a_k)))\\
    & = & {\rm PBW}_{\cB}(b\pi(a)\diff\pi(a_1).\diff\pi(a_2)\dots \diff\pi(a_k) )\\
    & = & {\rm PBW}_{\cB}(\pi_S(Y))\;.
\end{array}
$$
%
%
This proves the commutativity of the rightmost square. As for the leftmost square, we define $\theta$ by using the
restriction of the isomorphism $\Id_\cB\otimes\PBW_\cA$ to the ideal $(1_\cB\otimes\diff I)$: by commutativity of
the rightmost square the morphism $\Id_\cB\otimes\PBW_\cA$ sends $(1_\cB\otimes\diff I)$ in the image of
$\imath_U$. Consider the canonical surjection $\Gr_k'':I_\cB\cap(\cB\otimes U_k(\cA))\to I_\cB\cap(\cB\otimes
U_k(\cA)) / I_\cB\cap(\cB\otimes U_{k-1}(\cA))$. The graded morphism $\theta$ is given by
\begin{equation*}
 \theta(b\otimes a\;\diff a_1\diff a_2\dots\diff a_k)=\Gr_k''(b\otimes
 a\;\beta_\cA(a_1)\beta_\cA(a_2)\dots\beta_\cA(a_k))\;,
\end{equation*}%
where $b\in\cB$ and $a,a_1,\dots,a_k\in\cA$ are such that $a_i\in I$ for at least one index $i$. By construction,
the leftmost square is commutative.

\underline{Step 4:} \emph{$I_\cB$ as a left ideal.} We claim that $I_\cB$, which was defined as the ideal of
$\cB\otimes_\cA U(\cA)$ generated by $1_\cB\otimes\b_\cA(I)$ coincides with $I^L_\cB$, the \emph{left} ideal of
$\cB\otimes_\cA U(\cA)$ generated by $1_\cB\otimes\b_\cA(I)$. This property will be useful in the next step when we
prove that $\theta$ is surjective. Since $U(\cA)$ is generated by the images of $\a_\cA$ and $\b_\cA$ it suffices
to show that $I^L_\cB$ is stable for multiplication on the right by $1_\cB\otimes\a_\cA(a)$ and by
$1_\cB\otimes\b_\cA(a)$, for all $a\in \cA$. For $\jmath\in I$ we have
\begin{eqnarray*}
  (1_\cB\otimes\b_\cA(\jmath))\cdot(1_\cB\otimes\a_\cA(a))&=&
  1_\cB\otimes\b_\cA(\jmath).\a_\cA(a)\\
   & = & 1_{\cB}\otimes \alpha_{\cA}(a).\beta_{\cA}(j)+1_{\cB}\otimes \alpha_{\cA}(\{j,a\}) \\ 
 & = &(1_{\cB}\otimes \alpha_{\cA}(a))\cdot(1_{\cB}\otimes\beta_{\cA}(j))+\pi(\a_{\cA}(\{j,a\})\otimes1_{U_{\cA}}\\ 
  &=&(1_\cB\otimes\a_\cA(a))\cdot(1_\cB\otimes\b_\cA(\jmath))\;,
\end{eqnarray*}
which belongs to $I^L_\cB$; we have used in the last step that $\pb{\jmath,a}\in I$ (because $I$ is a Poisson
ideal). Similarly,
\begin{equation*}
  (1_\cB\otimes\b_\cA(\jmath))\cdot(1_\cB\otimes\b_\cA(a))=1_\cB\otimes\b_\cA\pb{\jmath,a}+
  (1_\cB\otimes\b_\cA(a))\cdot(1_\cB\otimes\b_\cA(\jmath))\in I^L_\cB\;.
\end{equation*}%
This proves our claim.

\underline{Step 5:} \emph{$\cB\otimes_\cA \Gr(U(\cA))$ has no non-trivial zero divisors.} Every $\cB$-module
becomes an $\cA$-module upon using $\pi:\cA\to \cB$ and similarly for every $\cB$-module morphism. Since $\cA$ is a
smooth algebra, $\Omega(\cA)$ is a projective $\cA$-module. It follows easily that $\cB\otimes_\cA\Omega(\cA)$ is a
projective $\cB$-module.  By definition, there exist $\cB$-modules $L$ and $N$, with $L$ free, such that
$L\simeq(\cB\otimes_\cA\Omega(\cA))\oplus N$. Since $L$ is free and since $\cB$ has no non-trivial zero divisors,
$\Sym_\cB L$ also has no non-trivial zero divisors. But $\Sym_\cB (\cB\otimes_\cA\Omega(\cA))$ is isomorphic to a
subalgebra of $\Sym_\cB L$, hence also has no non-trivial zero divisors. We can now conclude in view of the
isomorphisms (of $\cB$-modules)
\begin{equation*}
  \Sym_\cB(\cB\otimes_\cA \Omega(\cA))\simeq \cB\otimes_\cA \Sym_\cA\Omega(\cA)\simeq \cB\otimes_\cA \Gr(U(\cA))\;.
\end{equation*}%

\underline{Step 6:} \emph{Surjectivity of $\theta$.} For this step (only), we use our assumption that $I$ is
generated (as an ideal) by a single element, say $I=(\jmath)$, with $\jmath\in\cA$. Let $\Gr_k''(X')\in
I_\cB\cap(\cB\otimes_\cA U_k(\cA)) / I_\cB\cap(\cB\otimes U_{k-1}(\cA))$. If $X'\in I_\cB\cap(\cB\otimes
U_\ell(\cA))$ with $\ell<k$, then $\Gr_k''(X')=0=\theta(0)$. Thus we can suppose that $X'\in \cB\otimes_\cA
U_k(\cA)$ and $X'\notin \cB\otimes_\cA U_{k-1}(\cA)$.  We show that $\Gr_k''(X')=\theta(X)$ for some
$X\in(1_\cB\otimes\diff I)$ of degree $k$. Since $I_\cB$ is the left ideal generated by~$\jmath$ (see Step 4), we
can write $X'$ as $X'=Y'\cdot(1_\cB\otimes\b_\cA(\jmath))$, where $Y'\in \cB\otimes_\cA U(\cA)$. Since
$\cB\otimes_\cA \Gr(U(\cA))$ has no non-trivial zero divisors (Step~5), $Y'\in\cB\otimes_\cA U_{k-1}(\cA)$ and
$Y'\not\in \cB\otimes_\cA U_{k-2}(\cA)$. Therefore we can write $Y'$ as $Y'=Y_1'+Y'_2$ where $Y_1'\in\cB\otimes_\cA
U_{k-1}(\cA)$ is of the form $Y_1'=\sum_i b_i\otimes a_i\b_\cA(a_{1,i}).\b_\cA(a_{2,i})\dots \b_\cA(a_{k-1,i})$ and
$Y'_2\in \cB\otimes_\cA U_{k-2}(\cA)$. Then
\begin{eqnarray*}
  \Gr''_k(X')&=&\Gr''_k(Y'\cdot(1_\cB\otimes\b_\cA(\jmath)))=\Gr''_{k-1}(Y')\Gr''_1(1_\cB\otimes
  \b_\cA(\jmath))\\ &=&\Gr''_{k-1}(Y'_1)\Gr''_1(1_\cB\otimes\b_\cA(\jmath))\;.
\end{eqnarray*}%
Since $U(\cA)$ satisfies the PBW theorem, $Y'_1=\Id_\cB\otimes\PBW_\cA(Y)$ for some homogeneous element $Y\in
\cB\otimes_\cA\Sym_\cA(\Omega(\cA))$ of degree $k-1$. It follows that $\theta_k(Y\cdot (1_\cB\otimes\diff
\jmath))=\Gr_k''(Y'\cdot(1_\cB\otimes\b_\cA(\jmath)))$, so we can choose $X=Y\cdot (1_\cB\otimes\diff \jmath)$ to
obtain that $\Gr_k''(X')=\theta(X)$ with $X$ of degree $k$.
\end{proof}
In geometrical terms, the theorem is valid for arbitrary Poisson hypersurfaces of any smooth affine Poisson
variety. The example which follows is of this form.

\begin{example}\label{exa:Nambu_again}
As we have seen in Example \ref{exa:Nambu}, any pair of complex polynomials $(P,Q)$ defines a Poisson structure on
$\cA:=\bbC[X_1,X_2,X_3]$ by setting
\begin{equation}\label{eq:Nambu_end}
  \pb{X_1,X_2}:=Q\pp P{X_3}\;,\quad   \pb{X_2,X_3}:=Q\pp P{X_1}\;,\quad   \pb{X_3,X_1}:=Q\pp P{X_2}\;.
\end{equation}%
Since $P$ is a Casimir function of this Poisson structure, $(P)$ is a Poisson ideal of $\cA$ and
$\cB:=\bbC[X_1,X_2,X_3]/(P)$ is the algebra of functions of a Poisson surface, which may be a singular surface (for
example when $P$ is homogeneous of degree at least two). If $P$ is irreducible, so that $\cB$ is an integral
domain, then according to the above theorem, $\cB$ satisfies the PBW theorem. This example covers many classical
singular surfaces, such as the well-known Klein surfaces (see~\cite{AL}).
\end{example}
We also give an example of a Poisson algebra which does not satisfy the conditions of Theorem
\ref{thm:PBW_singular}, but yet the proof of this theorem can be used, with minor modifications, to show that it
satisfies the PBW theorem.

\begin{example}
We pick up again the previous example, but we take now for~$P$ the reducible polynomial $P:=X_1X_2X_3$. In this
case, $\cB:=\bbC[X_1,X_2,X_3]/(P)$ is the algebra of functions of a singular Poisson surface, which is the union of
the three coordinate planes in $\bbC^3$. Step 5 in the proof of Theorem \ref{thm:PBW_singular} is not valid
anymore, because $\cB$ now has non-trivial zero divisors. However, a close inspection of Step 6 in the proof
reveils that we only need to prove that
$1_\cB\otimes\diff\jmath\in\cB\otimes_\cA\Sym_\cA\Omega(\cA)\simeq\cB[Y_1,Y_2,Y_3]$ is not a zero divisor, i.e.,
that
\begin{equation*}
  \pi(X_2X_3)Y_1+\pi(X_1X_3)Y_2+\pi(X_1X_2)Y_3
\end{equation*}%
is not a zero divisor of $\cB[Y_1,Y_2,Y_3]$, where we recall that $\pi:\cA\to\cB$ denotes the canonical
surjection. Thus, we need to show that if $F$ is a polynomial in $Y_1,Y_2,Y_3$ with coefficients in $\cB$ and
\begin{equation}\label{eq:prod_zero}
  F(\pi(X_2X_3)Y_1+\pi(X_1X_3)Y_2+\pi(X_1X_2)Y_3)=0 \;,
\end{equation}  
then $F=0$. To do this, we may assume that $F$ is a homogenous polynomial of degree $n$ in the variables $Y_i$, say
$F=\sum_{i+j+k=n}\pi(F_{i,j,k})Y_1^{i}Y_2^{j}Y_3^{k}$. If (\ref{eq:prod_zero}) holds, then for all $n\geqs0$,
$$
  \sum_{i+j+k=n+1}\pi(F_{i-1,j,k}X_2X_3+F_{i,j-1,k}X_1X_3+ F_{i,j,k-1}X_1X_2)Y_1^{i}Y_2^{j}Y_3^{k}=0\;,
$$
hence for all $i,j,k\in\bbN$, the polynomial $X_1X_2X_3$ divides
$F_{i-1,j,k}X_2X_3+F_{i,j-1,k}X_1X_3+F_{i,j,k-1}X_1X_2$. Since $F_{i,j,k}=0$ whenever one of the indices $i,j,k$ is
negative, this means that $X_1$ divides every $F_{i,j,k}$. By symmetry, every $F_{i,j,k}$ is also divisible by
$X_2$ and by $X_3$. It follows that $\pi(F_{i,j,k})=0$ for all $i,j,k\geqs 0$, as was to be shown. We may now
conclude, in view of the proof of Theorem \ref{thm:PBW_singular} that $\cB$ satisfies the PBW theorem.
\end{example}

\subsection{The symmetrization map}
We now show that for a smooth Poisson algebra $\cA$, the PBW map $\Sym_\cA(\Omega(\cA))\to \Gr(U(\cA))$, which is
an isomorphism of $\cA$-algebras can be lifted to an isomorphism of $\cA$-modules $\Sym_\cA(\Omega(\cA))\to
U(\cA)$. To do this, we need to assume that we can divide elements of our base ring $R$ by arbitrary integers, so
we will assume in this paragraph that $R$ contains the field $\bbQ$ of rational numbers.

First, the following diagram is a commutative diagram of $\cA$-modules:
\begin{center}
  \begin{tikzcd}
    T^k_\cA(\Omega(\cA))\arrow{r}{\psi_k}\arrow{d}{\tau_k}&U_k(\cA)\arrow{d}{\Gr_k}\\
    \Sym_\cA^k(\Omega(\cA))\arrow{r}[swap]{\PBW_\cA^k}&\Gr^k(U(\cA))
  \end{tikzcd}
\end{center}
In this diagram, $\PBW_\cA^k$ is the restriction of $\PBW_\cA$ to $\Sym_\cA^k(\Omega(\cA))$, the homogenous
elements of degree $k$ of $\Sym_\cA(\Omega(\cA))$. The morphism $\tau_k$ is the canonical surjection and $\Gr_k$ is
the morphism which was introduced in (\ref{eq:rho}). Finally, the existence of a (non-unique) $\cA$-module morphism
$\psi_k$ making the diagram commutative follows from the fact that $T^k_\cA(\Omega(\cA))$ is a projective
$\cA$-module, which is in turn a consequence of the fact that $\cA$ is a smooth algebra (i.e., that $\Omega(\cA)$
is a projective $\cA$-module).
%

Let us denote by $T^{'k}_\cA(\Omega(\cA))\subset T^k_\cA(\Omega(\cA))$ the tensors which are symmetric, that is
invariant with respect to the standard action of the symmetric group~$S_k$; the restrictions of $\psi_k$ and
$\tau_k$ to this subspace are denoted by $\psi'_k$ and $\tau'_k$. Also, the image of $\psi'_k$ is denoted by
$U^k(\cA)$, because it can be viewed in a natural way as the degree $k$ component of a natural grading on
$U(\cA)$. Indeed, the above commutative diagram restricts to a new commutative diagram
\begin{center}
  \begin{tikzcd}
    T^{'k}_\cA(\Omega(\cA))\arrow{r}{\psi'_k}\arrow{d}{\tau'_k}&U^k(\cA)\arrow{d}{\Gr_k}\\
    \Sym_\cA^k(\Omega(\cA))\arrow{r}[swap]{\PBW_\cA^k}&\Gr^k(U(\cA))
  \end{tikzcd}
\end{center}
where $\tau'_k$ is now an isomorphism; since $\PBW_\cA^k$ is also an isomorphism, it follows that
$\Gr_k\circ\psi'_k: T^{'k}_\cA(\Omega(\cA))\to\Gr^k(U(\cA))$ is an isomorphism, and hence that the map $\psi'_k$
in this diagram is an isomorphism between $T^{'k}_\cA(\Omega(\cA))$ and a complement of $U_{k-1}(\cA)$ in
$U_k(\cA)$. As a corollary,
\begin{equation*}
  U_k(\cA)=U^k(\cA)\oplus U_{k-1}(\cA)=U^k(\cA)\oplus U^{k-1}(\cA)\oplus\cdots\oplus U^{0}(\cA)\;,
\end{equation*}%
and so $U(\cA)$ is graded by $\cA$-modules, $U(\cA)=\oplus_{i\in\bbN}U^i(\cA)$.

Since in the above diagram all maps are isomorphisms (of $\cA$-modules), we obtain by composition for every $k$ an
isomorphism of $\cA$-modules
\begin{equation*}
  \omega_k:\Sym_\cA^k(\Omega(\cA))\to U^k(\cA)\;,
\end{equation*}%
as claimed.
%


\bibliographystyle{abbrv} \bibliography{ref}

\end{document}